CrossMark

# Principal Lyapunov Exponents and Principal Floquet Spaces of Positive Random Dynamical Systems. III. Parabolic Equations and Delay Systems

**Janusz Mierczyński · Wenxian Shen**

**Abstract** This is the third part in a series of papers concerned with principal Lyapunov exponents and principal Floquet subspaces of positive random dynamical systems in ordered Banach spaces. The current part focuses on applications of general theory, developed in the authors' paper Mierczyński and Shen (Trans Am Math Soc 365(10):5329–5365, 2013), to positive continuous-time random dynamical systems on infinite dimensional ordered Banach spaces arising from random parabolic equations and random delay systems. It is shown under some quite general assumptions that measurable linear skew-product semidynamical systems generated by random parabolic equations and by cooperative systems of linear delay differential equations admit measurable families of generalized principal Floquet subspaces, and generalized principal Lyapunov exponents.

---

Dedicated to Professor John Mallet-Paret on the occasion of his 60th birthday.

---

J. Mierczyński (✉)
Institute of Mathematics and Computer Science, Wrocław University of Technology,
Wybrzeże Wyspiańskiego 27, 50-370 Wrocław, Poland
e-mail: mierczyn@pwr.edu.pl

W. Shen
Department of Mathematics and Statistics, Auburn University, Auburn, AL 36849, USA
e-mail: wenxish@mail.auburn.edu

## 1 Introduction

This is the third part of a series of several papers. The series is devoted to the study of principal Lyapunov exponents and principal Floquet subspaces of positive random dynamical systems in ordered Banach spaces.

Lyapunov exponents play an important role in the study of asymptotic dynamics of linear and nonlinear random evolution systems. Oseledets obtained in [22] important results on Lyapunov exponents and measurable invariant families of subspaces for finite-dimensional dynamical systems, which are called now the Oseledets multiplicative ergodic theorem. Since then a huge amount of research has been carried out toward alternative proofs of the Oseledets multiplicative ergodic theorem (see [2,10,11,16,20,23,24] and the references contained therein) and extensions of the Osedelets multiplicative theorem for finite dimensional systems to certain infinite dimensional ones (see [2,10–12,16,20,23–26], and references therein).

The largest finite Lyapunov exponents (or top Lyapunov exponents) and the associated invariant subspaces of both deterministic and random dynamical systems play special roles in the applications to nonlinear systems. Classically, the top finite Lyapunov exponent of a positive deterministic or random dynamical system in an ordered Banach space is called the *principal Lyapunov exponent* if the associated invariant family of subspaces corresponding to it consists of one-dimensional subspaces spanned by positive vectors (in such case, invariant subspaces are called the *principal Floquet subspaces*). For more on those subjects see [18].

In the first part of the series, [18], we introduced the notions of generalized principal Floquet subspaces, generalized principal Lyapunov exponents, and generalized exponential separations, which extend the corresponding classical notions. The classical theory of principal Lyapunov exponents, principal Floquet subspaces, and exponential separations for strongly positive and compact deterministic systems is extended to quite general positive random dynamical systems in ordered Banach spaces.

In the present, third part of the series, we consider applications of the general theory developed in [18] to positive random dynamical systems arising from random parabolic equations and systems of delay differential equations. To be more specific, let $((\Omega, \mathfrak{F}, \mathbb{P}), \theta_t)$ be an ergodic metric dynamical system. We consider a family, indexed by $\omega \in \Omega$, of second order partial differential equations

$$\frac{\partial u}{\partial t} = \sum_{i=1}^{N} \frac{\partial}{\partial x_i} \left( \sum_{j=1}^{N} a_{ij}(\theta_t \omega, x) \frac{\partial u}{\partial x_j} + a_i(\theta_t \omega, x) u \right) + \sum_{i=1}^{N} b_i(\theta_t \omega, x) \frac{\partial u}{\partial x_i}$$
$$+ c_0(\theta_t \omega, x) u, \qquad t > s, \ x \in D, \tag{1.1}$$

where $s \in \mathbb{R}$ is an initial time and $D \subset \mathbb{R}^N$ is a bounded domain with boundary $\partial D$, complemented with boundary condition

$$\mathcal{B}_{\theta_t \omega} u = 0, \quad t > s, \ x \in \partial D, \tag{1.2}$$

where

$$\mathcal{B}_\omega u = \begin{cases} u & \text{(Dirichlet)} \\ \sum_{i=1}^{N} \left( \sum_{j=1}^{N} a_{ij}(\omega, x) \frac{\partial u}{\partial x_j} + a_i(\omega, x) u \right) v_i + d_0(\omega, x) u & \text{(Robin)}. \end{cases}$$

Above, $\nu = (\nu_1, \ldots, \nu_N)$ denotes the unit normal vector pointing out of $\partial D$. We also consider the following system of linear random delay differential equations,

$$\frac{du}{dt} = A(\theta_t\omega)u(t) + B(\theta_t\omega)u(t-1) \tag{1.3}$$

where $u \in \mathbb{R}^N$, and $A(\omega)$, $B(\omega)$ are $N$ by $N$ real matrices (we write $A(\omega)$, $B(\omega) \in \mathbb{R}^{N \times N}$):

$$A(\omega) = \begin{pmatrix} a_{11}(\omega) & a_{12}(\omega) & \cdots & a_{1N}(\omega) \\ a_{21}(\omega) & a_{22}(\omega) & \cdots & a_{2N}(\omega) \\ \vdots & \vdots & \ddots & \vdots \\ a_{N1}(\omega) & a_{N2}(\omega) & \cdots & a_{NN}(\omega) \end{pmatrix},$$

and

$$B(\omega) = \begin{pmatrix} b_{11}(\omega) & b_{12}(\omega) & \cdots & b_{1N}(\omega) \\ b_{21}(\omega) & b_{22}(\omega) & \cdots & b_{2N}(\omega) \\ \vdots & \vdots & \ddots & \vdots \\ b_{N1}(\omega) & b_{N2}(\omega) & \cdots & b_{NN}(\omega) \end{pmatrix}.$$

Under quite general assumptions, it is shown that (1.1) + (1.2) generates a monotone random dynamical system on $X = L^2(D)$ and that (1.3) generates a monotone random dynamical system on $X = C([-1, 0], \mathbb{R}^N)$ provided that $A(\omega)$ and $B(\omega)$ are cooperative (see (OA1)).

Among others, we obtain the following results.

(1) Under some general assumptions, (1.1) + (1.2), (1.3) have nontrivial entire positive orbits (see Theorems 3.1, 4.1 for detail);

(2) Assume some focusing property. (1.1) + (1.2), (1.3) have measurable invariant families of one-dimensional subspaces $\{\tilde{E}_1(\omega)\}$ spanned by positive vectors (*generalized principal Floquet subspaces*) and whose associated Lyapunov exponent is the top Lyapunov exponent of the system (*generalized principal Lyapunov exponent*) (see Theorems 3.2, 4.2 for detail);

(3) Assume some strong focusing property. (1.1) + (1.2) have also measurable invariant families of one-codimensional subspaces which are exponentially separated from the generalized principal Floquet subspaces (see Theorem 3.3 for detail);

We remark that (1)–(3) are analogs of principal eigenvalue and principal eigenfunction theory for elliptic and periodic parabolic equations. Our main assumptions on (1.1) + (1.2) are the boundedness of $a_{ij}$, $a_i$, $b_i$ and $d_0$. No boundedness of $c_0$ is assumed. The results of the current paper hence extend the corresponding ones in [17] (it is assumed in [17] that $c_0$ is also bounded). In addition to the cooperative assumption, our main assumptions on (1.3) are the irreducibility of $B(\omega)$ or the positivity of $B(\omega)$. Such assumptions are also used in [21]. No boundedness of $A(\omega)$ and $B(\omega)$ is assumed in the current paper and the results of the current paper extend those in [21,27] for cooperative systems of delay differential equations.

It should be pointed out that the generalized principal Lyapunov exponents in (2) may be $-\infty$. In such a case, when generalized exponential separation holds, the (nontrivial) invariant measurable decomposition associated with the generalized exponential separation is essentially finer than the (trivial) decomposition in the Oseledets multiplicative ergodic theorem.

The results obtained in this paper would have important applications to the study of asymptotic dynamics of nonlinear random parabolic equations and systems of random delay differential equations.

The rest of this paper is organized as follows. First, for the reader's convenience, in Sect. 2 we recall some notions, assumptions, definitions, and main results established in part I [18]. We then consider random systems arising from parabolic equations and cooperative systems of delay differential equations in Sects. 3 and 4, respectively.

## 2 General Theory

In this section, we recall some general theory established in part I to be applied in this paper. To do so, we first introduce some notions, assumptions, and definitions introduced in part I. Then we recall some of the main results in part I.

2.1 Notions, Assumptions, and Definitions

In this subsection, we introduce some notions, assumptions, and definitions introduced in part I. The reader is referred to part I [18] for detail.

If $f$ is a real function defined on a set $Y$, we define its *nonnegative* (resp. *nonpositive*) *part* $f^+$ ($f^-$) as

$$f^+(y) := \begin{cases} f(y) & \text{if } f(y) \geq 0 \\ 0 & \text{if } f(y) < 0 \end{cases}, \qquad f^-(y) := \begin{cases} -f(y) & \text{if } f(y) \leq 0 \\ 0 & \text{if } f(y) > 0 \end{cases}.$$

For a metric space $(Y, d)$, $B(y; \epsilon)$ denotes the closed ball in $Y$ centered at $y \in Y$ and with radius $\epsilon > 0$. Further, $\mathfrak{B}(Y)$ stands for the $\sigma$-algebra of all Borel subsets of $Y$.

A *probability space* is a triple $(\Omega, \mathfrak{F}, \mathbb{P})$, where $\Omega$ is a set, $\mathfrak{F}$ is a $\sigma$-algebra of subsets of $\Omega$, and $\mathbb{P}$ is a probability measure defined for all $F \in \mathfrak{F}$. We always assume that the measure $\mathbb{P}$ is complete.

For a Banach space $X$, with norm $\|\cdot\|$, we will denote by $X^*$ its dual and by $\langle \cdot, \cdot \rangle$ the standard duality pairing (that is, for $u \in X$ and $u^* \in X^*$ the symbol $\langle u, u^* \rangle$ denotes the value of the bounded linear functional $u^*$ at $u$). Without further mention, we understand that the norm in $X^*$ is given by $\|u^*\| = \sup\{|\langle u, u^* \rangle| : \|u\| \leq 1\}$.

For Banach spaces $X_1$, $X_2$, $\mathcal{L}(X_1, X_2)$ stands for the Banach space of bounded linear mappings from $X_1$ into $X_2$, endowed with the standard norm. Instead of $\mathcal{L}(X, X)$ we write $\mathcal{L}(X)$.

$((\Omega, \mathfrak{F}, \mathbb{P}), (\theta_t)_{t \in \mathbb{R}})$ (we may simply write it as $(\theta_t)_{t \in \mathbb{R}}$, or as $(\theta_t)$) denotes an ergodic metric dynamical system.

For a metric dynamical system $((\Omega, \mathfrak{F}, \mathbb{P}), (\theta_t)_{t \in \mathbb{R}})$, $\Omega' \subset \Omega$ is *invariant* if $\theta_t(\Omega') = \Omega'$ for all $t \in \mathbb{R}$. $((\Omega, \mathfrak{F}, \mathbb{P}), (\theta_t)_{t \in \mathbb{R}})$ is said to be *ergodic* if for any invariant $F \in \mathfrak{F}$, either $\mathbb{P}(F) = 1$ or $\mathbb{P}(F) = 0$.

We write $\mathbb{R}^+$ for $[0, \infty)$. By a *measurable linear skew-product semidynamical system*, denoted by $\Phi = ((U_\omega(t))_{\omega \in \Omega, t \in \mathbb{R}^+}, (\theta_t)_{t \in \mathbb{R}})$, on a Banach space $X$ covering a metric dynamical system $(\theta_t)_{t \in \mathbb{R}}$ we understand a $(\mathfrak{B}(\mathbb{R}^+) \otimes \mathfrak{F} \otimes \mathfrak{B}(X), \mathfrak{B}(X))$-measurable mapping

$$\left[ \mathbb{R}^+ \times \Omega \times X \ni (t, \omega, u) \mapsto U_\omega(t)u \in X \right]$$

satisfying the following:

- 

$$U_\omega(0) = \text{Id}_X \quad \forall \omega \in \Omega, \tag{2.1}$$

$$U_{\theta_s \omega}(t) \circ U_\omega(s) = U_\omega(t + s) \qquad \forall \omega \in \Omega, \ t, s \in \mathbb{R}^+; \tag{2.2}$$

- for each $\omega \in \Omega$ and $t \in \mathbb{R}^+$, $[\, X \ni u \mapsto U_\omega(t)u \in X \,] \in \mathcal{L}(X)$.

Sometimes we write simply $\Phi = ((U_\omega(t)), (\theta_t))$.

Let $\Phi = ((U_\omega(t))_{\omega \in \Omega, t \in \mathbb{R}^+}, (\theta_t)_{t \in \mathbb{R}})$ be a measurable linear skew-product semidynamical system on a Banach space $X$ covering $(\theta_t)_{t \in \mathbb{R}}$. For $\omega \in \Omega$, $t \in \mathbb{R}^+$ and $u^* \in X^*$ we define $U_\omega^*(t)u^*$ by

$$\langle u, U_\omega^*(t)u^* \rangle = \langle U_{\theta_{-t}\omega}(t)u, u^* \rangle \qquad \text{for each } u \in X \tag{2.3}$$

(in other words, $U_\omega^*(t)$ is the mapping dual to $U_{\theta_{-t}\omega}(t)$).

In case where the mapping

$$\left[\, \mathbb{R}^+ \times \Omega \times X^* \ni (t, \omega, u^*) \mapsto U_\omega^*(t)u^* \in X^* \,\right]$$

is $(\mathfrak{B}(\mathbb{R}^+) \otimes \mathfrak{F} \otimes \mathfrak{B}(X^*), \mathfrak{B}(X^*))$-measurable, we will call the measurable linear skew-product semidynamical system $\Phi^* = ((U_\omega^*(t))_{\omega \in \Omega, t \in \mathbb{R}^+}, (\theta_{-t})_{t \in \mathbb{R}})$ on $X^*$ covering $(\theta_{-t})_{t \in \mathbb{R}}$ the *dual* of $\Phi$.

For measurable linear skew-product semiflows it often happens that for any $\omega \in \Omega$ and any $u \in X$ the mapping

$$\left[\, \mathbb{R}^+ \ni t \mapsto U_\omega(t)u \in X \,\right]$$

is continuous.

**Proposition 2.1** *Let $X$ be a separable Banach space. Assume that $(\theta_t)_{t \in \mathbb{R}}$ is a metric dynamical system and that*

$$\left[\, \mathbb{R}^+ \times \Omega \times X \ni (t, \omega, u) \mapsto U_\omega(t)u \in X \,\right]$$

*is a mapping satisfying* (2.1), (2.2), *such that the following holds:*

- *For any $t \in \mathbb{R}^+$ and $u \in X$ the mapping*

$$\left[\, \Omega \ni \omega \mapsto U_\omega(t)u \in X \,\right]$$

  *is $(\mathfrak{F}, \mathfrak{B}(X))$-measurable.*
- *For any $\omega \in \Omega$ and $t \in \mathbb{R}^+$, $[\, X \ni u \mapsto U_\omega(t)u \in X \,] \in \mathcal{L}(X)$.*
- *For any $\omega \in \Omega$ and $u \in X$ the mapping*

$$\left[\, \mathbb{R}^+ \ni t \mapsto U_\omega(t)u \in X \,\right]$$

  *is continuous.*

*Then $\Phi = ((U_\omega(t))_{\omega \in \Omega, t \in \mathbb{R}^+}, (\theta_t)_{t \in \mathbb{R}})$ is measurable linear skew-product semiflow having the following properties:*

(i) *For any $T > 0$ and any $u \in X$ the mapping*

$$\left[\, \Omega \ni \omega \mapsto \left[\, [0, T] \ni t \mapsto U_\omega(t)u \in X \,\right] \in C([0, T], X) \,\right]$$

  *is $(\mathfrak{F}, \mathfrak{B}(C([0, T], X)))$-measurable.*
(ii) *For any $\omega \in \Omega$ the mapping*

$$[\, \mathbb{R}^+ \times X \ni (t, u) \mapsto U_\omega(t)u \in X \,]$$

  *is continuous.*

*Proof* The fact that $[\, (t, \omega, u) \mapsto U_\omega(t)u \,]$ is $(\mathfrak{B}(\mathbb{R}^+) \otimes \mathfrak{F} \otimes \mathfrak{B}(X), \mathfrak{B}(X))$-measurable follows from [1, Lemma 4.51 on p. 153]. Part (i) is a consequence of [1, Theorem 4.55 on p. 155].

To prove (ii), fix $\omega \in \Omega$ and $T > 0$ and observe that for any $u \in X$ the set $\{\, \|U_\omega(t)u\| : t \in [0, T]\,\}$ is bounded. Hence, by the Uniform Boundedness Theorem, the set $\{\, \|U_\omega(t)\| : t \in [0, T]\,\}$ is bounded (by $M > 0$, say). Take a sequence $(t_n)_{n=1}^\infty \subset [0, T]$ convergent to $t$ and a sequence $(u_n)_{n=1}^\infty \subset X$ convergent to $u$. We estimate

$$\|U_\omega(t_n)u_n - U_\omega(t)u\| \leq \|U_\omega(t_n)(u_n - u)\| + \|U_\omega(t_n)u - U_\omega(t)u\|$$
$$\leq M\|u_n - u\| + \|U_\omega(t_n)u - U_\omega(t)u\|,$$

which goes to 0 as $n \to \infty$.                                                        $\square$

By a *cone* in a Banach space $X$ we understand a closed convex set $X^+$ such that

- $\alpha \geq 0$ and $u \in X^+$ imply $\alpha u \in X^+$, and
- $X^+ \cap (-X^+) = \{0\}$.

A pair $(X, X^+)$, where $X$ is a Banach space and $X^+$ is a cone in $X$, is referred to as an *ordered Banach space*.

If $(X, X^+)$ is an ordered Banach space, for $u, v \in X$ we write $u \leq v$ if $v - u \in X^+$, and $u < v$ if $u \leq v$ and $u \neq v$. The symbols $\geq$ and $>$ are used in an analogous way.

For a nonzero $u \in X^+$ we denote by $C_u$ the *component* of $u$: $v \in C_u$ if and only if there are positive numbers, $\underline{\alpha}, \overline{\alpha}$, such that $\underline{\alpha}v \leq u \leq \overline{\alpha}v$.

For an ordered Banach space $(X, X^+)$ denote by $(X^*)^+$ the set of all $u^* \in X^*$ such that $\langle u, u^* \rangle \geq 0$ for all $u \in X^+$. The set $(X^*)^+$ has the properties of a cone, except that $(X^*)^+ \cap (-(X^*)^+) = \{0\}$ need not be satisfied (such sets are called *wedges*).

If $(X^*)^+$ is a cone we call it the *dual cone*. This happens, for instance, when $X^+$ is total (that is, $X^+ + (-X^+)$ is dense in $X$).

Sometimes an ordered Banach space $(X, X^+)$ is a lattice: any two $u, v \in X$ have a least upper bound $u \vee v$ and a greatest lower bound $u \wedge v$. In such a case we write $u^+ := u \vee 0$, $u^- := (-u) \vee 0$, and $|u| := u^+ + u^-$. We have $u = u^+ - u^-$ for any $u \in X$.

An ordered Banach space $(X, X^+)$ being a lattice is a *Banach lattice* if there is a norm $\|\cdot\|$ on $X$ (a *lattice norm*) such that for any $u, v \in X$, if $|u| \leq |v|$ then $\|u\| \leq \|v\|$. From now on, when speaking of a Banach lattice we assume that the norm on $X$ is a lattice norm.

For application purposes, we give some examples of Banach lattices.

*Example 2.1* Let $D \subset \mathbb{R}^N$ be a bounded domain. The separable Banach space $L_p(D)$, where $1 < p < \infty$, is considered with the standard norm (denoted by $\|\cdot\|_p$). The standard cone $L_p(D)^+$ equals $\{\, u \in L_p(D) : u(x) \geq 0 \text{ for Lebesgue-a.e. } x \in D \,\}$. The pair $(L_p(D), L_p(D)^+)$ is a Banach lattice, and the norm $\|\cdot\|_p$ is a lattice norm. The dual cone in $L_p(D)^* = L_q(D)$, where $\frac{1}{p} + \frac{1}{q} = 1$, equals $L_q(D)^+$.

*Example 2.2* Let $N$ be a positive integer. Recall that $(\mathbb{R}^N, (\mathbb{R}^N)^+)$, where $(\mathbb{R}^N)^+$ denotes the set of vectors with nonnegative coordinates, is a Banach lattice, and that both the Euclidean norm $\|\cdot\|$ and the $\ell_1$-norm $\|\cdot\|_1$ on $\mathbb{R}^N$ are lattice norms (cf., e.g., [19, Example 2.1]).

We consider the separable Banach space $C([-1, 0], \mathbb{R}^N)$ of continuous $\mathbb{R}^N$-valued functions with norms

$$\|u\| = \sup_{t \in [-1, 0]} \|u(t)\|, \quad \|u\|_1 = \sup_{t \in [-1, 0]} \|u(t)\|_1 \quad \text{for} \quad u(\cdot) \in C\left([-1, 0], \mathbb{R}^N\right)$$

(the fact that, for instance, $\|\cdot\|$ is used for both the Euclidean norm on $\mathbb{R}^N$ and the corresponding supremum norm on $C([-1, 0], \mathbb{R}^N)$ should not cause any misunderstanding). We define

$$C([-1, 0], \mathbb{R}^N)^+ := \left\{ u \in C\left([-1, 0], \mathbb{R}^N\right) : u(\tau) \in (\mathbb{R}^N)^+ \text{ for all } \tau \in [-1, 0] \right\}.$$

$\left(C([-1,0],\mathbb{R}^N), C\left([-1,0],\mathbb{R}^N\right)^+\right)$ is a Banach lattice: for $u \in C\left([-1,0],\mathbb{R}^N\right)$, $u^+$ (resp. $u^-$) are defined by $u^+(\tau) := (u(\tau))^+$ (resp. $u^-(\tau) := (u(\tau))^-$) for all $\tau \in [-1,0]$. Both norms $\|\cdot\|$ and $\|\cdot\|_1$ are lattice norms.

Note that the dual Banach space $C([-1,0],\mathbb{R}^N)^*$ is not separable.

We introduce now our assumptions.

**(C0)** (Banach lattice) $(X, X^+)$ *is a separable Banach lattice with* dim $X \geq 2$.

Observe that if, $(X, X^+)$ satisfies (C0), then all (A0)(i), (A0)(ii), and (A0)(iii) in [18] are satisfied.

**(C0)\*** (Banach lattice) $(X^*, (X^*)^+)$ *is a separable Banach lattice with* dim $X^* \geq 2$.

**(C1)** (Integrability/injectivity/complete continuity) $\Phi = (U_\omega(t), (\theta_t))$ *is a measurable linear skew-product semidynamical system* on a separable Banach space $X$ *covering an ergodic metric dynamical system* $(\theta_t)$ *on* $(\Omega, \mathfrak{F}, \mathbb{P})$, *with the complete measure* $\mathbb{P}$, *satisfying the following*:

(i) (Integrability) *The functions*

$$\left[\Omega \ni \omega \mapsto \sup_{0 \leq s \leq 1} \ln^+ \|U_\omega(s)\| \in [0, \infty)\right] \in L_1((\Omega, \mathfrak{F}, \mathbb{P}))$$

and

$$\left[\Omega \ni \omega \mapsto \sup_{0 \leq s \leq 1} \ln^+ \|U_{\theta_s\omega}(1-s)\| \in [0, \infty)\right] \in L_1((\Omega, \mathfrak{F}, \mathbb{P})).$$

(ii) (Injectivity) *For each* $\omega \in \Omega$ *the linear operator* $U_\omega(1)$ *is injective.*

(iii) (Complete continuity) *For each* $\omega \in \Omega$ *the linear operator* $U_\omega(1)$ *is completely continuous.*

In the sequel, by (C1)\*(i), (C1)\*(ii) and (C1)\*(iii) we will understand the counterparts of (C1)(i), (C1)(ii) and (C1)(iii) for the dual measurable linear skew-product semidynamical system $\Phi^*$. More precisely, for example (C1)\*(ii) means the following: "the mapping $\left[\mathbb{R}^+ \times \Omega \times X^* \ni (t, \omega, u^*) \mapsto U_\omega^*(t)u^* \in X^*\right]$ is $\left(\mathcal{B}(\mathbb{R}^+) \otimes \mathfrak{F} \otimes \mathcal{B}(X^*), \mathcal{B}(X^*)\right)$-measurable, and for each $\omega \in \Omega$ the linear operator $U_\omega^*(1)$ is injective."

Observe that, assuming that the measurability in the definition of $\Phi^*$ holds, if (C1)(i) is satisfied then (C1)\*(i) is satisfied, too; similarly, if (C1)(iii) is satisfied then (C1)(iii)\* is satisfied.

**(C2)** (Positivity) $(X, X^+)$ *satisfies* (C0) *and* $\Phi = ((U_\omega(t)), (\theta_t))$ *is a measurable linear skew-product semidynamical system on* $X$ *covering an ergodic metric dynamical system* $(\theta_t)$ *on* $(\Omega, \mathfrak{F}, \mathbb{P})$, *satisfying the following*:

$$U_\omega(t)u_1 \leq U_\omega(t)u_2$$

*for any* $\omega \in \Omega$, $t \geq 0$ *and* $u_1, u_2 \in X$ *with* $u_1 \leq u_2$.

**(C2)\*** (Positivity) $(X^*, (X^*)^+)$ *satisfies* (C0)\* *and* $\Phi^* = \left((U_\omega^*(t)), (\theta_{-t})\right)$ *is a* measurable linear skew-product semidynamical system *on* $X^*$ *covering an ergodic metric dynamical system* $(\theta_{-t})$ *on* $(\Omega, \mathfrak{F}, \mathbb{P})$, *satisfying the following*:

$$U_\omega^*(t)u_1^* \leq U_\omega^*(t)u_2^*$$

*for any* $\omega \in \Omega$, $t \geq 0$ *and* $u_1^*, u_2^* \in X^*$ *with* $u_1^* \leq u_2^*$.

**(C3)** (Focusing) (C2) *is satisfied and there are* $\mathbf{e} \in X^+$ *with* $\|\mathbf{e}\| = 1$ *and an* $(\mathfrak{F}, \mathcal{B}(\mathbb{R}))$-*measurable function* $\varkappa \colon \Omega \to [1, \infty)$ *with* $\ln^+ \ln \varkappa \in L_1((\Omega, \mathfrak{F}, \mathbb{P}))$ *such that for any* $\omega \in \Omega$ *and any nonzero $u \in X^+$ there is $\beta(\omega, u) > 0$ with the property that*

$$\beta(\omega, u)\mathbf{e} \le U_\omega(1)u \le \varkappa(\omega)\beta(\omega, u)\mathbf{e}.$$

**(C3)\*** *(Focusing) (C2)\* is satisfied and there are* $\mathbf{e}^* \in (X^*)^+$ *with* $\|\mathbf{e}^*\| = 1$ *and an* $(\mathfrak{F}, \mathfrak{B}(\mathbb{R}))$*-measurable function* $\varkappa^* \colon \Omega \to [1, \infty)$ *with* $\ln^+ \ln \varkappa^* \in L_1((\Omega, \mathfrak{F}, \mathbb{P}))$ *such that for any* $\omega \in \Omega$ *and any nonzero* $u^* \in (X^*)^+$ *there is* $\beta^*(\omega, u^*) > 0$ *with the property that*

$$\beta^* \left(\omega, u^*\right) \mathbf{e}^* \le U_\omega^*(1)u^* \le \varkappa^*(\omega)\beta^* \left(\omega, u^*\right) \mathbf{e}^*.$$

**(C4)** *(Strong focusing) (C3), (C3)\* are satisfied and* $\ln \varkappa \in L_1((\Omega, \mathfrak{F}, \mathbb{P}))$, $\ln \varkappa^* \in L_1((\Omega, \mathfrak{F}, \mathbb{P}))$, *and* $\langle \mathbf{e}, \mathbf{e}^* \rangle > 0$.

**(C5)** *(Strong positivity in one direction) There are* $\bar{\mathbf{e}} \in X^+$ *with* $\|\bar{\mathbf{e}}\| = 1$ *and an* $(\mathfrak{F}, \mathfrak{B}(\mathbb{R}))$*-measurable function* $\nu \colon \Omega \to (0, \infty)$, *with* $\ln^- \nu \in L_1((\Omega, \mathfrak{F}, \mathbb{P}))$, *such that*

$$U_\omega(1)\bar{\mathbf{e}} \ge \nu(\omega)\bar{\mathbf{e}} \quad \forall\, \omega \in \Omega.$$

**(C5)\*** *(Strong positivity in one direction) There are* $\bar{\mathbf{e}}^* \in (X^*)^+$ *with* $\|\bar{\mathbf{e}}^*\| = 1$ *and an* $(\mathfrak{F}, \mathfrak{B}(\mathbb{R}))$*-measurable function* $\nu^* \colon \Omega \to (0, \infty)$, *with* $\ln^- \nu^* \in L_1((\Omega, \mathfrak{F}, \mathbb{P}))$, *such that*

$$U_\omega^*(1)\bar{\mathbf{e}}^* \ge \nu^*(\omega)\bar{\mathbf{e}}^* \quad \forall\, \omega \in \Omega.$$

*Remark 2.1* We can replace time 1 with some $T > 0$ in (C1), (C3), (C4), (C5), and (C1)\*, (C3)\*, (C5)\*.

We now state the definitions introduced in [18]. Throughout the rest of this subsection, until revocation, we assume (C0) and (C2).

**Definition 2.1** (*Entire positive orbit*) For $\omega \in \Omega$, by an *entire positive orbit* of $U_\omega$ we understand a mapping $v_\omega \colon \mathbb{R} \to X^+$ such that $v_\omega(s + t) = U_{\theta_s\omega}(t)v_\omega(s)$ for any $s \in \mathbb{R}$ and $t \in \mathbb{R}^+$. The function constantly equal to zero is referred to as the *trivial entire orbit*.

Entire positive orbits of $\Phi^*$ are defined in a similar way.

A family $\{E(\omega)\}_{\omega \in \Omega_0}$ of $l$-dimensional vector subspaces of $X$ is *measurable* if there are $(\mathfrak{F}, \mathfrak{B}(X))$-measurable functions $v_1, \ldots, v_l \colon \Omega_0 \to X$ such that $(v_1(\omega), \ldots, v_l(\omega))$ forms a basis of $E(\omega)$ for each $\omega \in \Omega_0$.

Let $\{E(\omega)\}_{\omega \in \Omega_0}$ be a family of $l$-dimensional vector subspaces of $X$, and let $\{F(\omega)\}_{\omega \in \Omega_0}$ be a family of $l$-codimensional closed vector subspaces of $X$, such that $E(\omega) \oplus F(\omega) = X$ for all $\omega \in \Omega_0$. We define the *family of projections associated with the decomposition* $E(\omega) \oplus F(\omega) = X$ as $\{P(\omega)\}_{\omega \in \Omega_0}$, where $P(\omega)$ is the linear projection of $X$ onto $F(\omega)$ along $E(\omega)$, for each $\omega \in \Omega_0$.

The family of projections associated with the decomposition $E(\omega) \oplus F(\omega) = X$ is called *strongly measurable* if for each $u \in X$ the mapping $[\, \Omega_0 \ni \omega \mapsto P(\omega)u \in X \,]$ is $(\mathfrak{F}, \mathfrak{B}(X))$-measurable.

We say that the decomposition $E(\omega) \oplus F(\omega) = X$, with $\{E(\omega)\}_{\omega \in \Omega_0}$ finite-dimensional, is *invariant* if $\Omega_0$ is invariant, $U_\omega(t)E(\omega) = E(\theta_t\omega)$ and $U_\omega(t)F(\omega) \subset F(\theta_t\omega)$, for each $t \in \mathbb{T}^+$.

A strongly measurable family of projections associated with the invariant decomposition $E(\omega) \oplus F(\omega) = X$ is referred to as *tempered* if

$$\lim_{t \to \pm\infty} \frac{\ln \|P(\theta_t\omega)\|}{t} = 0 \quad \mathbb{P}\text{-a.e. on } \Omega_0.$$

**Definition 2.2** (*Generalized principal Floquet subspaces and principal Lyapunov exponent*)
A family of one-dimensional subspaces $\{\tilde{E}(\omega)\}_{\omega \in \tilde{\Omega}}$ of $X$ is called a family of *generalized principal Floquet subspaces* of $\Phi = ((U_\omega(t)), (\theta_t))$ if $\tilde{\Omega} \subset \Omega$ is invariant, $\mathbb{P}(\tilde{\Omega}) = 1$, and

(i)  $\tilde{E}(\omega) = \text{span}\,\{w(\omega)\}$ with $w \colon \tilde{\Omega} \to X^+ \backslash \{0\}$ being $(\mathfrak{F}, \mathfrak{B}(X))$-measurable,
(ii)  $U_\omega(t)\tilde{E}(\omega) = \tilde{E}(\theta_t\omega)$, for any $\omega \in \tilde{\Omega}$ and any $t > 0$,
(iii)  there is $\tilde{\lambda} \in [-\infty, \infty)$ such that

$$\tilde{\lambda} = \lim_{t \to \infty} \frac{\ln \|U_\omega(t)w(\omega)\|}{t} \quad \forall \omega \in \tilde{\Omega},$$

(iv)

$$\limsup_{t \to \infty} \frac{\ln \|U_\omega(t)u\|}{t} \leq \tilde{\lambda} \quad \forall \omega \in \tilde{\Omega},\ \forall u \in X \backslash \{0\}.$$

$\tilde{\lambda}$ is called the *generalized principal Lyapunov exponent* of $\Phi$ associated to the generalized principal Floquet subspaces $\{\tilde{E}(\omega)\}_{\omega \in \tilde{\Omega}}$.

Observe that if $\{\tilde{E}(\omega)\}_{\tilde{\omega} \in \Omega}$ is a family of generalized principal Floquet subspaces of $((U_\omega(t))_{\omega \in \Omega, t \in \mathbb{R}^+}, (\theta_t)_{t \in \mathbb{R}})$, then for any $\omega \in \tilde{\Omega}$, $v_\omega(\cdot)$ is an entire positive orbit, where

$$v_\omega(t) = \begin{cases} U_\omega(t)w(\omega), & t \geq 0 \\ \dfrac{\|w(\omega)\|}{\|U_{\theta_t\omega}(-t)w(\theta_t\omega)\|} w(\theta_t\omega), & t < 0. \end{cases}$$

In the literature on random linear skew-product dynamical systems the concept of the *top* (or the *largest*) *Lyapunov exponent* is introduced. It can be defined either as the largest exponential growth rate of the norms of the *individual* vectors (in such a case, when $\Phi$ has a family of generalized principal Floquet subspaces then the generalized principal Lyapunov exponent is, by definition, the top Lyapunov exponent), or as the exponential growth rate of the norms of the *operators*. These definitions are equivalent, however we have been unable to locate a concise proof in the existing literature. This is the reason why we decided to formulate and prove the result below (the proof is patterned after the proof of [13, Theorem 2.2], in the light of the first and second remarks on p. 528 of [13]).

**Proposition 2.2** *Assume that* $\Phi = ((U_\omega(t)), (\theta_t))$ *has a family of generalized principal Floquet subspaces, with the generalized principal Lyapunov exponent* $\tilde{\lambda}$. *Assume moreover (C1)(i). Then*

$$\lim_{t \to \infty} \frac{\ln \|U_\omega(t)\|}{t} = \tilde{\lambda}$$

*for any* $\omega \in \tilde{\Omega}$, *where* $\tilde{\Omega}$ *is as in the Definition 2.2.*

*Proof* We start by proving that

$$\lim_{n \to \infty} \frac{\ln \|U_\omega(n)\|}{n} = \tilde{\lambda} \tag{2.4}$$

for any $\omega \in \tilde{\Omega}$.

Fix some $\omega \in \tilde{\Omega}$ and $\lambda > \tilde{\lambda}$, and define functions $p_n \colon X \to [0, \infty)$, $n = 1, 2, \ldots$, and $p \colon X \to [0, \infty)$ as

$$p_n(u) := \frac{\|U_\omega(n)u\|}{e^{\lambda n}}, \quad p(u) := \sup_{n \in \mathbb{N}} p_n(u).$$

For $m = 1, 2, \ldots$ put

$$W_m := \{\, u \in X : p(u) \le m \,\} = \bigcap_{n=1}^{\infty} \{\, u \in X : p_n(u) \le m \,\}.$$

The sets $W_m$ are closed and their union equals the whole of $X$ (by Definition 2.2(iv)). By the Baire theorem, there is $m_0 \in \mathbb{N}$ such that $W_{m_0}$ has nonempty interior. In other words, there exist $v \in X$ and $\epsilon > 0$ such that $B(v; \epsilon) \subset W_{m_0}$. From this it follows that

$$\|U_\omega(n)(v + w)\| \le m_0 e^{\lambda n}$$

for all $w \in X$ with $\|w\| \le \epsilon$ and all $n = 1, 2, \ldots$. In particular, by taking $w = 0$ we have that $\|U_\omega(n)v\| \le m_0 e^{\lambda n}$ for all $n$. By the triangle inequality,

$$\|U_\omega(n)w\| \le 2m_0 e^{\lambda n}$$

for all $w \in X$ with $\|w\| \le \epsilon$ and all $n = 1, 2, \ldots$.

As $\lambda > \tilde{\lambda}$ is arbitrary, we have that

$$\limsup_{n \to \infty} \frac{\ln \|U_\omega(n)\|}{n} \le \tilde{\lambda},$$

which, combined with Definition 2.2(iii), gives (2.4).

The passage to the continuous time, under (C1)(i), goes by a standard argument, as presented for instance in the proof of [12, Lemma 3.4]. □

**Definition 2.3** (*Generalized exponential separation*) $\Phi = ((U_\omega(t)), (\theta_t))$ admits *a generalized exponential separation* if it has a family of generalized principal Floquet subspaces $\{\tilde{E}(\omega)\}_{\omega \in \tilde{\Omega}}$ and a family of one-codimensional subspaces $\{\tilde{F}(\omega)\}_{\omega \in \tilde{\Omega}}$ of $X$ satisfying the following

(i) $\tilde{F}(\omega) \cap X^+ = \{0\}$ for any $\omega \in \tilde{\Omega}$,
(ii) $X = \tilde{E}(\omega) \oplus \tilde{F}(\omega)$ for any $\omega \in \tilde{\Omega}$, where the decomposition is invariant, and the family of projections associated with this decomposition is strongly measurable and tempered,
(iii) there exists $\tilde{\sigma} \in (0, \infty]$ such that

$$\lim_{t \to \infty} \frac{1}{t} \ln \frac{\|U_\omega(t)|_{\tilde{F}(\omega)}\|}{\|U_\omega(t)w(\omega)\|} = -\tilde{\sigma} \quad \forall \omega \in \tilde{\Omega}.$$

We say that $\{\tilde{E}(\cdot), \tilde{F}(\cdot), \tilde{\sigma}\}$ *generates a generalized exponential separation*.

We end this subsection with the following proposition which follows from the Oseledets-type theorems proved in [12] (we do not assume (C0) or (C2) now).

**Proposition 2.3** *Let $X$ be a separable Banach space of infinite dimension. Let $\Phi = ((U_\omega(t)), (\theta_t))$ be a measurable linear skew-product semidynamical system satisfying (C1)(i)–(iii). Then there exist: an invariant $\Omega_0 \subset \Omega$, $\mathbb{P}(\Omega_0) = 1$, and $\lambda_1 \in [-\infty, \infty)$ such that*

$$\lim_{t \to \infty} \frac{\ln \|U_\omega(t)\|}{t} = \lambda_1 \quad \forall \omega \in \Omega_0. \tag{2.5}$$

*Moreover, if $\lambda_1 > -\infty$ then there are a measurable family $\{E_1(\omega)\}_{\omega \in \Omega_0}$ of vector subspaces of finite dimension, and a family $\{\hat{F}_1(\omega)\}_{\omega \in \Omega_0}$ of closed vector subspaces of finite codimension such that*

(i) *$X = E_1(\omega) \oplus \hat{F}_1(\omega)$ for any $\omega \in \Omega_0$, where the decomposition is invariant, and the family of projections associated with this decomposition is strongly measurable and tempered,*

(ii)
$$\lim_{t\to\pm\infty}\frac{\ln\|U_\omega(t)|_{E_1(\omega)}\|}{t}=\lim_{t\to\pm\infty}\frac{\ln\|U_\omega(t)u\|}{t}=\lambda_1\quad\forall\omega\in\Omega_0,\ u\in E_1(\omega)\backslash\{0\},$$

(iii)
$$\lim_{t\to\infty}\frac{\ln\|U_\omega(t)u\|}{t}=\lambda_1\quad\forall\omega\in\Omega_0,\ u\in E_1(\omega)\backslash\hat{F}_1(\omega),$$

(iv) *there is $\hat{\lambda}_2\in[-\infty,\lambda_1)$ such that*
$$\lim_{t\to\infty}\frac{\ln\|U_\omega(t)|_{\hat{F}_1(\omega)}\|}{t}=\hat{\lambda}_2\quad\forall\omega\in\Omega_0.$$

*Proof* See [18, Theorem 3.4, and (3.1) on p. 5342].                                            $\square$

2.2 General Theorems

In this subsection, we state some general theorems, most of which are established in part I.

The first theorem is on the existence of entire positive orbits.

**Theorem 2.1** (Entire positive orbits) *Assume* (C0), (C1)(i)–(iii) *and* (C2). *If $\lambda_1>-\infty(\lambda_1$ is as in (2.5)), then there is a measurable set $\Omega_1$ of $\omega\in\Omega_0$ with $\mathbb{P}(\Omega_1)=1$ such that for each $\omega\in\Omega_1$ there exists a nontrivial entire positive orbit $v_\omega\colon\mathbb{R}\to X^+$ of $\Phi=((U_\omega(t)),(\theta_t))$ and*
$$\lim_{t\to\pm\infty}\frac{\ln\|U_\omega(t)v_\omega\|}{t}=\lambda_1.$$

The above theorem follows from [18, Theorem 3.5].

Next theorem shows the existence of generalized Floquet subspaces and principal Lyapunov exponent.

**Theorem 2.2** (Generalized principal Floquet subspace and Lyapunov exponent) *Assume* (C0), (C1)(i), (C2) *and* (C3). *Then there exist an invariant set $\tilde{\Omega}_1\subset\Omega$, $\mathbb{P}(\tilde{\Omega}_1)=1$, and an $(\mathfrak{F},\mathfrak{B}(X))$-measurable function $w\colon\tilde{\Omega}_1\to X$, $w(\omega)\in C_{\mathbf{e}}$ and $\|w(\omega)\|=1$ for all $\omega\in\tilde{\Omega}_1$, having the following properties:*

(1)
$$w(\theta_t\omega)=\frac{U_\omega(t)w(\omega)}{\|U_\omega(t)w(\omega)\|}$$

*for any $\omega\in\tilde{\Omega}_1$ and $t\geq0$.*

(2) *Let for some $\omega\in\tilde{\Omega}_1$ a function $v_\omega\colon\mathbb{R}\to X^+\backslash\{0\}$ be a nontrivial entire positive orbit of $U_\omega$. Then $v_\omega(t)=\|v_\omega(0)\|w_\omega(t)$ for all $t\in\mathbb{R}$, where*
$$w_\omega(t):=\begin{cases}\left(U_{\theta_t\omega}(-t)|_{\tilde{E}_1(\theta_t\omega)}\right)^{-1}w(\omega)&\text{for }t<0\\U_\omega(t)w(\omega)&\text{for }t\geq0,\end{cases}$$

*with $\tilde{E}_1(\omega)=\mathrm{span}\{w(\omega)\}$.*

(3) *There exists $\tilde{\lambda}_1\in[-\infty,\infty)$ such that*
$$\tilde{\lambda}_1=\lim_{t\to\pm\infty}\frac{\ln\rho_t(\omega)}{t}=\int_\Omega\ln\rho_1\,d\mathbb{P}$$

*for each $\omega \in \tilde{\Omega}_1$, where*

$$\rho_t(\omega) := \begin{cases} \|U_\omega(t)w(\omega)\| & \text{for } t \geq 0, \\ 1/\|U_{\theta_t\omega}(-t)w(\theta_t\omega)\| & \text{for } t < 0. \end{cases}$$

(4) *For any $\omega \in \tilde{\Omega}_1$ and any $u \in X^+ \backslash \{0\}$ there holds*

$$\lim_{t \to \infty} \frac{\ln \|U_\omega(t)u\|}{t} = \tilde{\lambda}_1.$$

(5) *For $\omega \in \tilde{\Omega}_1$ and any $u \in X \backslash \{0\}$,*

$$\limsup_{t \to \infty} \frac{\ln \|U_\omega(t)u\|}{t} \leq \tilde{\lambda}_1,$$

*and then $\{\tilde{E}_1(\omega)\}_{\omega \in \tilde{\Omega}_1}$ is a family of generalized Floquet subspaces, with the generalized principal Lyapunov exponent equal to $\tilde{\lambda}_1$, and for any $\omega \in \tilde{\Omega}_1$ there holds*

$$\lim_{t \to \infty} \frac{\ln \|U_\omega(t)\|}{t} = \tilde{\lambda}_1.$$

(6) *Assume, moreover, that (C1)(ii)–(iii) hold. Then $\lambda_1$ in Proposition 2.3 equals $\tilde{\lambda}_1$.*
(7) *Assume, moreover, that (C5) holds. Then $\tilde{\lambda}_1 > -\infty$.*

*Proof* First of all, parts (1) through (3) are reformulations of [18, Theorem 3.6(1)–(3)].

We next prove (4). By (C3) and part (1), for each $\omega \in \tilde{\Omega}_1$ there are $\gamma_1(\omega) > 0$ and $\gamma_2(\omega) > 0$ such that

$$\gamma_1(\omega)\mathbf{e} \leq w(\omega) \leq \gamma_2(\omega)\mathbf{e},$$

which gives, via the monotonicity of the norm $\|\cdot\|$, that

$$\lim_{t \to \infty} \frac{\ln \|U_\omega(t)\mathbf{e}\|}{t} = \lim_{t \to \infty} \frac{\ln \|U_\omega(t)w(\omega)\|}{t} = \tilde{\lambda}_1.$$

Further, for each $\omega \in \Omega$ and each $u \in X^+ \backslash \{0\}$ there are $\tilde{\gamma}_1(\omega, u) > 0$ and $\tilde{\gamma}_2(\omega, u) > 0$ such that

$$\tilde{\gamma}_1(\omega, u)\mathbf{e} \leq U_\omega(1)u \leq \tilde{\gamma}_2(\omega, u)\mathbf{e},$$

which again gives that

$$\lim_{t \to \infty} \frac{\ln \|U_\omega(t)u\|}{t} = \lim_{t \to \infty} \frac{\ln \|U_{\theta_1\omega}(t-1)U_\omega(1)u\|}{t} = \lim_{t \to \infty} \frac{\ln \|U_\omega(t)\mathbf{e}\|}{t} = \tilde{\lambda}_1.$$

Now we prove (5). By (C0) and (C2), for any $\omega \in \tilde{\Omega}_1$, $t > 0$ and $u \in X$,

$$\|U_\omega(t)u\| \leq \|U_\omega(t)|u|\|.$$

It then follows from (4) that

$$\limsup_{t \to \infty} \frac{\ln \|U_\omega(t)u\|}{t} \leq \tilde{\lambda}_1,$$

and then $\{\tilde{E}_1(\omega)\}_{\omega \in \tilde{\Omega}_1}$ is a family of generalized Floquet subspaces, with the generalized principal Lyapunov exponent equal to $\tilde{\lambda}_1$. Hence by Proposition 2.2, for any $\omega \in \tilde{\Omega}_1$ there holds

$$\lim_{t \to \infty} \frac{\ln \|U_\omega(t)\|}{t} = \tilde{\lambda}_1.$$

Finally, (6) follows from [18, Theorem 3.6(4)] and (7) follows from the arguments of [18, Theorem 3.8(6)].                                                                                         □

The theorem below is a counterpart of Theorem 2.2 for the dual system.

**Theorem 2.3** (Generalized principal Floquet subspace and Lyapunov exponent) *Assume* (C0)*, (C1)*(i), (C2)* *and* (C3)*. *Then there exist an invariant set* $\tilde{\Omega}_1^* \subset \Omega$, $\mathbb{P}(\tilde{\Omega}_1^*) = 1$, *and an* $(\mathfrak{F}, \mathfrak{B}(X^*))$-*measurable function* $w^* \colon \tilde{\Omega}_1^* \to X^*$, $w^*(\omega) \in C_{\mathbf{e}^*}$ *and* $\|w^*(\omega)\| = 1$ *for all* $\omega \in \tilde{\Omega}_1^*$, *having the following properties:*

(1)
$$w^*(\theta_{-t}\omega) = \frac{U_\omega^*(t)w^*(\omega)}{\|U_\omega^*(t)w^*(\omega)\|}$$

*for any* $\omega \in \tilde{\Omega}_1^*$ *and* $t \geq 0$.

(2) *Let for some* $\omega \in \tilde{\Omega}_1^*$ *a function* $v_\omega^* \colon \mathbb{R} \to (X^*)^+ \backslash \{0\}$ *be a nontrivial entire positive orbit of* $U_\omega^*$. *Then* $v_\omega^*(t) = \|v_\omega^*(0)\| w_\omega^*(t)$ *for all* $t \in \mathbb{R}$, *where*

$$w_\omega^*(t) := \begin{cases} \left(U_{\theta_{-t}\omega}^*(-t)|_{\tilde{E}_1^*(\theta_{-t}\omega)}\right)^{-1} w^*(\omega) & \text{for } t < 0 \\ U_\omega^*(t)w^*(\omega) & \text{for } t \geq 0, \end{cases}$$

*where* $\tilde{E}_1^*(\omega) = \mathrm{span}\{w^*(\omega)\}$.

(3) *There exists* $\tilde{\lambda}_1^* \in [-\infty, \infty)$ *such that*

$$\tilde{\lambda}_1^* = \lim_{t \to \pm\infty} \frac{\ln \rho_t^*(\omega)}{t} = \int_\Omega \ln \rho_1^* \, d\mathbb{P}$$

*for each* $\omega \in \tilde{\Omega}_1$, *where*

$$\rho_t^*(\omega) := \begin{cases} \|U_\omega^*(t)w^*(\omega)\| & \text{for } t \geq 0, \\ 1/\|U_{\theta_{-t}\omega}^*(-t)w^*(\theta_{-t}\omega)\| & \text{for } t < 0, \end{cases}$$

(4) *If* (C0), (C1)(i), (C2) *and* (C3) *are satisfied, then* $\tilde{\lambda}_1 = \tilde{\lambda}_1^*$.

*Proof* It is just a restatement of [18, Theorem 3.7].                                   □

To state the next theorem, we introduce some notions. Assume (C0), (C1)(i), (C2), (C3), and (C0)*, (C1)*(i), (C2)* and (C3)*. Let $\tilde{\Omega}_1$, $w(\cdot)$, and $\tilde{\Omega}_1^*$, $w^*(\cdot)$ be as in Theorems 2.2 and 2.3, respectively. For $\omega \in \tilde{\Omega}_1^*$, define $\tilde{F}_1(\omega) := \{u \in X : \langle u, w^*(\omega) \rangle = 0\}$. Then $\{\tilde{F}_1(\omega)\}_{\omega \in \tilde{\Omega}_1^*}$ is a family of one-codimensional subspaces of $X$, such that $U_\omega(t)\tilde{F}_1(\omega) \subset \tilde{F}_1(\theta_t\omega)$ for any $\omega \in \tilde{\Omega}_1^*$ and any $t \geq 0$.

**Theorem 2.4** (Generalized exponential separation) *Assume* (C0), (C1)(i), (C2), (C0)*, (C1)*(i), (C2)*, *and* (C4). *Then there is an invariant set* $\tilde{\Omega}_0$, $\mathbb{P}(\tilde{\Omega}_0) = 1$, *having the following properties.*

(1) *The family* $\{\tilde{P}(\omega)\}_{\omega \in \tilde{\Omega}_0}$ *of projections associated with the invariant decomposition* $\tilde{E}_1(\omega) \oplus \tilde{F}_1(\omega) = X$ *is strongly measurable and tempered.*

(2) $\tilde{F}_1(\omega) \cap X^+ = \{0\}$ *for any* $\omega \in \tilde{\Omega}_0$.

(3) *For any $\omega \in \tilde{\Omega}_0$ and any $u \in X \setminus \tilde{F}_1(\omega)$ (in particular, for any nonzero $u \in X^+$) there holds*

$$\lim_{t \to \infty} \frac{\ln \|U_\omega(t)\|}{t} = \lim_{t \to \infty} \frac{\ln \|U_\omega(t)u\|}{t} = \tilde{\lambda}_1.$$

(4) *There exist $\tilde{\sigma} \in (0, \infty]$ and $\tilde{\lambda}_2 \in [-\infty, \infty)$, $\tilde{\lambda}_2 = \tilde{\lambda}_1 - \tilde{\sigma}$, such that*

$$\lim_{t \to \infty} \frac{1}{t} \ln \frac{\|U_\omega(t)|_{\tilde{F}(\omega)}\|}{\|U_\omega(t)w(\omega)\|} = -\tilde{\sigma}$$

*and*

$$\lim_{t \to \infty} \frac{\ln \|U_\omega(t)|_{\tilde{F}_1(\omega)}\|}{t} = \tilde{\lambda}_2$$

*for each $\omega \in \tilde{\Omega}_0$. Hence $\Phi$ admits a generalized exponential separation.*

(5) *Assume moreover* (C1)(ii)–(iii) *and* (C1)*(ii)–(iii). *If $\tilde{\lambda}_1 > -\infty$ then, in the notation of Proposition 2.3, $\tilde{\lambda}_2 = \hat{\lambda}_2 (< \hat{\lambda}_1)$ and $E_1(\omega) = \hat{E}_1(\omega)$ and $\hat{F}_1(\omega) = \tilde{F}_1(\omega)$ for $\mathbb{P}$-a.e. $\omega \in \Omega_0$.*

*Proof* It follows from [18, Theorem 3.8]. □

## 3 Linear Random Parabolic Equations

In this section, we consider applications of the general results stated in Sect. 2 to linear random parabolic equations.

Let $((\Omega, \mathfrak{F}, \mathbb{P}), (\theta_t)_{t \in \mathbb{R}})$ be an ergodic metric dynamical system, with $\mathbb{P}$ complete. Consider (1.1)+(1.2), that is, a family, indexed by $\omega \in \Omega$, of second order partial differential equations,

$$\frac{\partial u}{\partial t} = \sum_{i=1}^N \frac{\partial}{\partial x_i} \left( \sum_{j=1}^N a_{ij}(\theta_t \omega, x) \frac{\partial u}{\partial x_j} + a_i(\theta_t \omega, x) u \right) + \sum_{i=1}^N b_i(\theta_t \omega, x) \frac{\partial u}{\partial x_i}$$
$$+ c_0(\theta_t \omega, x) u, \qquad t > s, \ x \in D, \tag{3.1}$$

where $s \in \mathbb{R}$ is an initial time and $D \subset \mathbb{R}^N$ is a bounded domain with boundary $\partial D$, complemented with boundary condition

$$\mathcal{B}_{\theta_t \omega} u = 0, \quad t > s, \ x \in \partial D, \tag{3.2}$$

where

$$\mathcal{B}_\omega u = \begin{cases} u & \text{(Dirichlet)} \\ \displaystyle\sum_{i=1}^N \left( \sum_{j=1}^N a_{ij}(\omega, x) \frac{\partial u}{\partial x_j} + a_i(\omega, x) u \right) \nu_i + d_0(\omega, x) u & \text{(Robin).} \end{cases}$$

Above, $\nu = (\nu_1, \ldots, \nu_N)$ denotes the unit normal vector pointing out of $\partial D$. When $d_0 \equiv 0$ in the Robin case, $\mathcal{B}_\omega u = 0$ is also referred to as the Neumann boundary condition.

In addition, we consider also the *adjoint problem* to (3.1) + (3.2), that is,

$$-\frac{\partial u}{\partial s} = \sum_{i=1}^N \frac{\partial}{\partial x_i} \left( \sum_{j=1}^N a_{ji}(\theta_s \omega, x) \frac{\partial u}{\partial x_j} - b_i(\theta_s \omega, x) u \right) - \sum_{i=1}^N a_i(\theta_s \omega, x) \frac{\partial u}{\partial x_i}$$
$$+ c_0(\theta_s \omega, x) u, \qquad s < t, \ x \in D, \tag{3.3}$$

where $t \in \mathbb{R}$ is a final time, complemented with boundary condition

$$\mathcal{B}^*_{\theta_s \omega} u = 0, \quad s < t, \ x \in \partial D, \tag{3.4}$$

where $\mathcal{B}^*_\omega = \mathcal{B}_\omega$ in the Dirichlet boundary conditions case, or

$$\mathcal{B}^*_\omega u = \sum_{i=1}^{N} \left( \sum_{j=1}^{N} a_{ji}(\omega, x) \frac{\partial u}{\partial x_j} - b_i(\omega, x) u \right) v_i + d_0(\omega, x) u \tag{3.5}$$

in the Robin case.

When we want to emphasize that (3.1) + (3.2) is considered for some (fixed) $\omega \in \Omega$ we write $(3.1)_\omega$ + $(3.2)_\omega$. The same holds for (3.3) + (3.4).

Throughout the present section, $\|\cdot\|$ stands for the norm in $L_2(D)$ or for the norm in $\mathcal{L}(L_2(D))$, depending on the context. Sometimes we use summation convention. For example, we can write (3.1) as

$$\partial_t u = \partial_i \big( a_{ij}(\theta_t \omega, x) \partial_j u + a_i(\theta_t \omega, x) u \big) + b_i(\theta_t \omega, x) \partial_i u + c_0(\theta_t \omega, x) u.$$

When speaking of properties satisfied by points in $D$, we use the expression "for a.e. $x \in D$" to indicate that the $N$-dimensional Lebesgue measure of the set of points not satisfying the property is zero. Similarly, when speaking of properties satisfied by points in $\partial D$, we use the expression "for a.e. $x \in \partial D$" to indicate that the $(N-1)$-dimensional Lebesgue measure of the set of points not satisfying the property is zero. The expressions "for a.e. $(t, x) \in \mathbb{R} \times D$," "for a.e. $(t, x) \in \mathbb{R} \times \partial D$" are used in an analogous way.

## 3.1 Measurable Linear Skew-Product Semiflows

In this subsection, we give a sketch of the existence theory for (weak) $L_2(D)$-solutions of (3.1) + (3.2) (or of (3.3) + (3.4)). It is an appropriate modification of the proof presented in the authors' monograph [17, Chapter 2 and Subsection 4.1.1].

First of all, we introduce some assumptions on $D$ and the coefficients of the problem (3.1) + (3.2).

**(PA0)** (Boundary regularity) $D \subset \mathbb{R}^N$ *is a bounded domain with Lipschitz boundary $\partial D$.*

**(PA1)** (Measurability) *The functions $a_{ij} \colon \Omega \times D \to \mathbb{R}$ $(i, j = 1, \ldots, N)$, $a_i \colon \Omega \times D \to \mathbb{R}$ $(i = 1, \ldots, N)$, $b_i \colon \Omega \times D \to \mathbb{R}$ $(i = 1, \ldots, N)$ and $c_0 \colon \Omega \times D \to \mathbb{R}$ are $(\mathfrak{F} \otimes \mathfrak{B}(D), \mathfrak{B}(\mathbb{R}))$-measurable. In the case of Robin boundary conditions the function $d_0 \colon \Omega \times \partial D \to [0, \infty)$ is $(\mathfrak{F} \otimes \mathfrak{B}(\partial D), \mathfrak{B}(\mathbb{R}))$-measurable.*

**(PA2)**

(i) (Boundedness of second and first order terms) *For each $\omega \in \Omega$ the functions $[ \, (t, x) \mapsto a_{ij}(\theta_t \omega, x) \, ]$ $(i, j = 1, \ldots N)$, $[ \, (t, x) \mapsto a_i(\theta_t \omega, x) \, ]$ $(i = 1, \ldots, N)$ and $[ \, (t, x) \mapsto b_i(\theta_t \omega, x) \, ]$ $(i = 1, \ldots, N)$ belong to $L_\infty(\mathbb{R} \times D)$, with their $L_\infty(\mathbb{R} \times D)$-norms bounded uniformly in $\omega \in \Omega$. In the Robin case, for each $\omega \in \Omega$ the functions $[ \, (t, x) \mapsto d_0(\theta_t \omega, x) \, ]$ belong to $L_\infty(\mathbb{R} \times \partial D)$, with their $L_\infty(\mathbb{R} \times \partial D)$-norms bounded uniformly in $\omega \in \Omega$.*

(ii) (Local boundedness of zero order terms) *For each $\omega \in \Omega$, $c_0(\omega, \cdot) \in L_\infty(D)$. Moreover, there are $(\mathfrak{F}, \mathfrak{B}(\mathbb{R}))$-measurable functions $c_0^{(+)} \colon \Omega \to [0, \infty)$, $c_0^{(-)} \colon \Omega \to (-\infty, 0]$, such that for each $\omega \in \Omega$,*

- $$c_0^{(-)}(\omega) \le c_0(\omega, x) \le c_0^{(+)}(\omega), \qquad x \in D,$$

- *the mappings* $\left[\, \mathbb{R} \ni t \mapsto c_0^{(\pm)}(\theta_t\omega) \in \mathbb{R} \,\right]$ *are continuous.*

**(PA3)** (Ellipticity) *There exists* $\alpha_0 > 0$ *such that for each* $\omega \in \Omega$ *there holds*

$$\sum_{i,j=1}^{N} a_{ij}(\theta_t\omega, x)\xi_i\xi_j \geq \alpha_0 \sum_{i=1}^{N} \xi_i^2, \qquad \xi = (\xi_1, \ldots, \xi_N) \in \mathbb{R}^N,$$

*and*

$$a_{ij}(\theta_t\omega, x) = a_{ji}(\theta_t\omega, x), \qquad i, j = 1, 2, \ldots, N,$$

*for a.e.* $(t, x) \in \mathbb{R} \times D$.

For $\omega \in \Omega$ define functions $a_{ij}^{\omega} : \mathbb{R} \times D \to \mathbb{R}, i, j = 1, \ldots, N$, by $a_{ij}^{\omega}(t, x) := a_{ij}(\theta_t\omega, x)$, and similarly for $a_i^{\omega}, b_i^{\omega}(i = 1, \ldots, N), c_0, d_0$. Put

$$a^{\omega} := \left( (a_{ij}^{\omega})_{i,j=1}^{N}, (a_i^{\omega})_{i=1}^{N}, (b_i^{\omega})_{i=1}^{N}, c_0^{\omega}, d_0^{\omega} \right)$$

(in the Dirichlet case we put $d_0$ constantly equal to zero).

For $s \in \mathbb{R}$, $M > 0$ and $T > 0$, let

$$\Omega_{s,M,T} := \left\{ \omega \in \Omega : -M \leq c_0^{(-)}(\theta_t\omega) \leq c_0^{(+)}(\theta_t\omega) \leq M \text{ for } s \leq t \leq s + T \right\}$$

**Lemma 3.1** *Assume* (PA0), (PA1) *and* (PA2)(ii). *Then for any* $s \in \mathbb{R}$, *any* $M > 0$ *and any* $T > 0$, *the set* $\Omega_{s,M,T}$ *is a measurable subset of* $\Omega$.

*Proof* Let $\mathbb{Q}$ be the set of all rational numbers. Then

$$\Omega_{s,M,T} = \bigcap_{\substack{s \leq r \leq s+T \\ r \in \mathbb{Q}}} \left\{ \omega : -M \leq c_0^{(-)}(\theta_r\omega) \leq c_0^{(+)}(\theta_r\omega) \leq M \right\}.$$

Clearly, for any $r \in \mathbb{Q}$, $\{\omega : -M \leq c_0^{(-)}(\theta_r\omega) \leq c_0^{(+)}(\theta_r\omega) \leq M\}$ is a measurable subset of $\Omega$. It then follows that $\Omega_{s,M,T}$ is measurable. $\qquad\square$

Observe that for any given $s \in \mathbb{R}$ and $T > 0$,

$$\Omega = \bigcup_{\substack{M > 0 \\ M \in \mathbb{Q}}} \Omega_{s,M,T}. \tag{3.6}$$

From now on we assume that (PA0) through (PA3) are satisfied.

For any $s \in \mathbb{R}$, $T > 0$ and $\omega \in \Omega$, the restriction

$$a_{s,T}^{\omega} := \left( (a_{ij}^{\omega})_{i,j=1}^{N}|_{(s,s+T) \times D}, (a_i^{\omega})_{i=1}^{N}|_{(s,s+T) \times D}, (b_i^{\omega})_{i=1}^{N}|_{(s,s+T) \times D}, c_0^{\omega}|_{(s,s+T) \times D}, \right.$$
$$\left. d_0^{\omega}|_{(s,s+T) \times \partial D} \right)$$

belongs to $L_{\infty}((s, s + T) \times D, \mathbb{R}^{N^2+2N+1}) \times L_{\infty}((s, s + T) \times \partial D, \mathbb{R})$. Moreover, the set $\{a_{s,T}^{\omega} : \omega \in \Omega_{s,M,T}\}$ is bounded in $L_{\infty}((s, s+T) \times D, \mathbb{R}^{N^2+2N+1}) \times L_{\infty}((s, s+T) \times \partial D, \mathbb{R})$.

For any $s \in \mathbb{R}$, $M > 0$ and $T > 0$ we write $Y_{s,M,T}$ for the closure of $\{a_{s,T}^{\omega} : \omega \in \Omega_{s,M,T}\}$ in the weak-* topology of $L_{\infty}((s, s + T) \times D, \mathbb{R}^{N^2+2N+1}) \times L_{\infty}((s, s + T) \times \partial D, \mathbb{R})$. $Y_{s,M,T}$ is a compact metrizable space.

We consider, for each $\tilde{a} = \left( (\tilde{a}_{ij})_{i,j=1}^{N}, (\tilde{a}_i)_{i=1}^{N}, (\tilde{b}_i)_{i=1}^{N}, \tilde{c}_0, \tilde{d}_0 \right) \in Y_{s,M,T}$,

$$\partial_t u = \partial_i \left( \tilde{a}_{ij}(t, x) \partial_j u + \tilde{a}_i(t, x) u \right) + \tilde{b}_i(t, x) \partial_i u + \tilde{c}_0(t, x) u, \ t \in (s, s + T], \ x \in D \tag{3.7}$$

complemented with boundary conditions

$$\mathcal{B}_{\tilde{a}}u = 0, \quad t \in (s, s + T], \ x \in \partial D, \tag{3.8}$$

where $\mathcal{B}_{\tilde{a}}u = u$ in the Dirichlet case, or

$$\mathcal{B}_{\tilde{a}}u = \big(\tilde{a}_{ij}(t, x)\partial_j u + \tilde{a}_i(t, x)u\big)v_i + \tilde{d}_0(t, x)u \tag{3.9}$$

in the Robin case. Recall that, if $\tilde{d}_0 \equiv 0$ in the Robin case, $\mathcal{B}_{\tilde{a}}u = 0$ is also referred to as the Neumann boundary condition. To emphasize the dependence of the equation on the parameter $\tilde{a}$ we write $(3.7)_{\tilde{a}} + (3.8)_{\tilde{a}}$.

Let $V$ be defined as follows

$$V := \begin{cases} \overset{\circ}{W}{}^1_2(D) & \text{(Dirichlet)} \\ W^1_2(D) & \text{(Neumann)} \\ W^1_{2,2}(D, \partial D) & \text{(Robin)} \end{cases} \tag{3.10}$$

where $\overset{\circ}{W}{}^1_2(D)$ is the closure of $\mathcal{D}(D)$ in $W^1_2(D)$ and $W^1_{2,2}(D, \partial D)$ is the completion of

$$V_0 := \big\{ v \in W^1_2(D) \cap C(\bar{D}) \cap C^\infty(D) : \|v\|_V < \infty \big\}$$

with respect to the norm $\|v\|_V := (\|\nabla v\|_2^2 + \|v\|_{2,\partial D}^2)^{1/2}$.

Let

$$W = W(s, s + T; V, V^*) := \{ v \in L_2((s, s + T), V) : \dot{v} \in L_2((s, s + T), V^*) \} \tag{3.11}$$

equipped with the norm

$$\|v\|_W := \left( \int_s^{s+T} \|v(\tau)\|_V^2 \, d\tau + \int_s^{s+T} \|\dot{v}(\tau)\|_{V^*}^2 \, d\tau \right)^{\frac{1}{2}},$$

where $\dot{v} := dv/dt$ is the time derivative in the sense of distributions taking values in $V^*$ (see [5, Chapter XVIII] for definitions).

For $\tilde{a} \in Y_{s,M,T}$ denote by $B_{\tilde{a}} = B_{\tilde{a}}(t, \cdot, \cdot)$ the bilinear form on $V$ associated with $\tilde{a}$,

$$B_{\tilde{a}}(t, u, v) := \int_D \big((\tilde{a}_{ij}(t, x)\partial_j u + \tilde{a}_i(t, x)u)\partial_i v - (\tilde{b}_i(t, x)\partial_i u + \tilde{c}_0(t, x)u)v\big) \, dx, \quad u, v \in V, \tag{3.12}$$

in the Dirichlet and Neumann boundary condition cases, and

$$B_{\tilde{a}}(t, u, v) := \int_D \big((\tilde{a}_{ij}(t, x)\partial_j u + \tilde{a}_i(t, x)u)\partial_i v - (\tilde{b}_i(t, x)\partial_i u + \tilde{c}_0(t, x)u)v\big) \, dx$$
$$+ \int_{\partial D} \tilde{d}_0(t, x)uv \, dH_{N-1}, \quad u, v \in V, \tag{3.13}$$

in the Robin boundary condition case, where $H_{N-1}$ stands for $(N-1)$-dimensional Hausdorff measure (since $\partial D$ is Lipschitz, $H_{N-1}$ is in fact Lebesgue surface measure).

**Definition 3.1** (*Weak solution*) Let $\tilde{a} \in Y_{s,M,T}$ and $u_0 \in L_2(D)$. A function $u \in L_2((s, s + T), V)$ is a *weak solution* of $(3.7)_{\tilde{a}} + (3.8)_{\tilde{a}}$ on $[s, s + T]$ with initial condition $u(s) = u_0$ if

$$- \int_s^{s+T} \langle u(\tau), v \rangle \dot{\phi}(\tau) \, d\tau + \int_s^{s+T} B_{\tilde{a}}(\tau, u(\tau), v) \phi(\tau) \, d\tau - \langle u_0, v \rangle \phi(s) = 0 \qquad (3.14)$$

for all $v \in V$ and $\phi \in \mathcal{D}([s, s + T))$, where $\mathcal{D}([s, s + T))$ is the space of all smooth real functions having compact support in $[s, s + T)$.

Our next assumptions will guarantee continuous dependence of solutions on parameters.

**(PA4)** (*Convergence almost everywhere*) *For any* $s \in \mathbb{R}$, $M > 0$ *and* $T > 0$, *for any sequence* $(\tilde{a}^{(n)})$ *converging in* $Y_{s,M,T}$ *to* $\tilde{a}$ *we have that* $\tilde{a}_{ij}^{(n)} \to \tilde{a}_{ij}$, $\tilde{a}_i^{(n)} \to \tilde{a}_i$, $\tilde{b}_i^{(n)} \to \tilde{b}_i$ *pointwise a.e. on* $(s, s + T) \times D$ *and (in the Robin boundary condition case)* $\tilde{d}_0^{(n)} \to \tilde{d}_0$ *pointwise a.e. on* $(s, s + T) \times \partial D$.

Below, in Sect. 3.2 we give sufficient conditions ((R)) for (PA4) to be fulfilled.

**Proposition 3.1** *Let* $M > 0$, $s \in \mathbb{R}$ *and* $T > 0$. *Then for each* $\tilde{a} \in Y_{s,M,T}$ *and each* $u_0 \in L_2(D)$ *there exists a unique weak solution* $\tilde{u}^{M,T}(\cdot; s, \tilde{a}, u_0)$ *of* $(3.7)_{\tilde{a}} + (3.8)_{\tilde{a}}$ *on* $[s, s + T]$ *with initial condition* $\tilde{u}(0) = u_0$, *satisfying the following properties.*

(i) (*Continuity in* $L_2(D)$) *For each* $\tilde{a} \in Y_{s,M,T}$ *the mapping*
$$\left[ [s, s + T] \times L_2(D) \ni (t, u_0) \mapsto \tilde{u}^{M,T}(t; s, \tilde{a}, u_0) \in L_2(D) \right]$$
*is continuous.*

(ii) (*Continuity in* $Y_{s,M,T}$) *Assume moreover* (PA4). *For each* $u_0 \in L_2(D)$ *the mapping*
$$\left[ (s, s + T] \times Y_{s,M,T} \ni (t, \tilde{a}) \mapsto \tilde{u}^{M,T}(t; s, \tilde{a}, u_0) \in L_2(D) \right]$$
*is continuous.*

(iii) (*Positivity*) *For any* $\tilde{a} \in Y_{s,M,T}$ *and any nonzero* $u_0 \in L_2(D)^+$ *there holds* $\tilde{u}^{M,T}(t; s, \tilde{a}, u_0) \in L_2(D)^+ \backslash \{0\}$ *for all* $t \in (s, s + T]$.

(iv) (*Compactness*) *Assume moreover* (PA4). *For any* $t_1, t_2 \in (s, s + T]$, $t_1 \leq t_2$, *and any bounded subset* $E \subset L_2(D)$, $\{ \tilde{u}^{M,T}(t; s, \tilde{a}, u_0) : \tilde{a} \in Y_{s,M,T}, \ t \in [t_1, t_2], \ u_0 \in E \}$ *is a relatively compact subset of* $L_2(D)$.

*Indication of proof.* Existence of weak solutions follows from [4, Theorem 2.4]. (i) follows from [4, Theorem 5.1]. (It should be remarked that, formally, in [4] the coefficients are assumed to be essentially bounded on $(s, \infty) \times D$ (resp. on $(s, \infty) \times D$). However, in the results utilized here only the values on $(s, s + T) \times D$ (resp. on $(s, s + T) \times D$) are of importance.)

(ii), (iii) and (iv) follow along the lines of the proof of [17, Theorem 2.4.1(3), Proposition 2.2.9, Proposition 2.2.5], respectively.                                                          □

**Proposition 3.2** *Let* $M > 0$, $s \in \mathbb{R}$ *and* $T > 0$. *Then for any two* $\tilde{a}^{(1)}, \tilde{a}^{(2)} \in Y_{s,M,T}$ *such that* $\tilde{a}_{ij}^{(1)}(t, x) = \tilde{a}_{ij}^{(2)}(t, x)$, $\tilde{a}_i^{(1)}(t, x) = \tilde{a}_i^{(2)}(t, x)$, $\tilde{b}_i^{(1)}(t, x) = \tilde{b}_i^{(2)}(t, x)$ *a.e. on* $(s, s + T) \times D$, *(in the Robin boundary condition case)* $\tilde{d}_0^{(1)}(t, x) = \tilde{d}_0^{(2)}(t, x)$ *a.e. on* $(s, s + T) \times \partial D$, *but* $\tilde{c}_0^{(1)}(t, x) \leq \tilde{c}_0^{(2)}(t, x)$ *a.e. on* $(s, S + T) \times D$, *and any* $u_0 \in L_2(D)^+$ *there holds*

$$\tilde{u}^{M,T}(t; s, \tilde{a}^{(1)}, u_0) \leq \tilde{u}^{M,T}(t; s, \tilde{a}^{(2)}, u_0)$$

*for all* $t \in (s, s + T]$, *where* $\leq$ *is understood in the* $L_2(D)^+$-*sense.*

*Proof* Compare the proof of [17, Proposition 2.2.10(1)–(2)].                    □

Define $\mathcal{E}_{s,M,T} : \Omega_{s,M,T} \to Y_{s,M,T}$ by $\mathcal{E}_{s,M,T}(\omega) := a^\omega$, restricted to $(s, s + T) \times D$ (or to $(s, s + T) \times \partial D$).

**Lemma 3.2** *Assume moreover* (PA4). *The mapping* $\mathcal{E}_{s,M,T}$ *is* $\big(\mathfrak{F}, \mathfrak{B}(Y_{s,M,T})\big)$-*measurable.*

*Proof* Cf. [17, Lemma 4.1.1].                                                   □

Fix for the moment $\omega \in \Omega$, $s \in \mathbb{R}$ and $u_0 \in L_2(D)$. We proceed now to the definition of the global weak solution $u(\cdot; s, \omega, u_0)$ of $(3.1)_\omega + (3.2)_\omega$ satisfying the initial condition $u(s) = u_0$.

For $n = 1, 2, 3, \ldots$ denote

$$M_n := \max\big\{ -\inf\{ c_0^{(-)}(\theta_t\omega) : t \in [s, s + n] \}, \sup\{ c_0^{(+)}(\theta_t\omega) : t \in [s, s + n] \} \big\}.$$

By (PA2)(ii), $M_n < \infty$.

**Definition 3.2** (*Global weak solution*) Let $\omega \in \Omega$, $s \in \mathbb{R}$ and $u_0 \in L_2(D)$. A *global weak solution* of $(3.1)_\omega + (3.2)_\omega$ with initial condition $u(s) = u_0$ is defined as

$$u(t; s, \omega, u_0) := \tilde{u}^{M_{\lfloor t-s \rfloor +1}, \lfloor t-s \rfloor +1}\big( t; s, \mathcal{E}_{s, M_{\lfloor t-s \rfloor +1}, \lfloor t-s \rfloor +1}(\omega), u_0 \big)$$

for $t > s$.

It follows from the uniqueness of weak solutions that the $u(\cdot; s, \omega, u_0)$ is well defined.

**Lemma 3.3** *The global weak solutions of* (3.1) + (3.2) *have the following properties.*

(i) (Time translation) *For any* $\omega \in \Omega$, *any* $s \le t$ *and any* $u_0 \in L_2(D)$ *there holds*

$$u(t; s, \omega, u_0) = u(t - s; 0, \theta_s\omega, u_0).$$

(ii) (Cocycle identity) *For any* $\omega \in \Omega$, *any* $s \le t_1 \le t_2$ *and any* $u_0 \in L_2(D)$ *there holds*

$$u(t_2; s, \omega, u_0) = u(t_2 - t_1; t_1, \theta_{t_1-s}\omega, u(t_1; s, \omega, u_0)). \tag{3.15}$$

*Indication of proof.* The proof goes by appropriately rewriting the proofs of Propositions 2.1.6 through 2.1.8 in [17].                                           □

Similarly, for $\omega \in \Omega$, $t \in \mathbb{R}$ and $u_0 \in L_2(D)$ we define a global weak solution $u^*(\cdot; t, \omega, u_0^*)$ of the adjoint problem (3.3) + (3.4) satisfying the final condition $u(t, \cdot) = u_0$.

**Lemma 3.4** *Let* $\omega \in \Omega$ *and* $s < t$. *Then for any* $u_0, u_0^* \in L_2(D)$ *there holds*

$$\langle u(t; s, \omega, u_0), u_0^* \rangle = \langle u_0, u^*(s; t, \omega, u_0^*) \rangle.$$

*Proof* See [17, Proposition 2.3.2].                                             □

We define

$$U_\omega(t)u_0 := u(t, 0; \omega, u_0), \qquad (\omega \in \Omega, \ t \ge 0, \ u_0 \in L_2(D)), \tag{3.16}$$
$$U_\omega^*(t)u_0^* := u^*(-t, 0; \omega, u_0^*), \qquad \big(\omega \in \Omega, \ t \ge 0, \ u_0^* \in L_2(D)\big). \tag{3.17}$$

From now on, we assume additionally that (PA4) is satisfied.

**Proposition 3.3** $\Phi = ((U_\omega(t))_{\omega \in \Omega, t \in \mathbb{R}^+}, (\theta_t)_{t \in \mathbb{R}})$ *is a measurable linear skew-product semidynamical system on* $L_2(D)$ *covering a metric dynamical system* $(\theta_t)_{t \in \mathbb{R}}$, *with* $\Phi^* = \left( (U_\omega^*(t))_{\omega \in \Omega, t \in \mathbb{R}^+}, (\theta_{-t})_{t \in \mathbb{R}} \right)$ *being its dual.*

*Proof* Equation (2.1) follows in a straightforward way from the definition of a weak solution, and (2.2) is a consequence of Lemma 3.3. The property that $U_\omega(t)$ belongs to $\mathcal{L}(L_2(D))$ follows from Proposition 3.1(i).

By arguments similar to those in the proof of [17, Lemma 4.1.3], for fixed $M > 0$ and $T > 0$ the mapping

$$\left[ [0, T] \times \Omega_{0,M,T} \times L_2(D) \ni (t, \omega, u_0) \mapsto u(t, \cdot; \omega, u_0) \in L_2(D) \right]$$

is $\left( \mathfrak{B}([0, T]) \otimes \mathfrak{F}|_{\Omega_{0,M,T}} \otimes \mathfrak{B}(L_2(D)), \mathfrak{B}(L_2(D)) \right)$-measurable. As $T > 0$ is arbitrary, it follows via (3.6) that the mapping

$$\left[ \mathbb{R}^+ \times \Omega \times L_2(D) \ni (t, \omega, u_0) \mapsto U_\omega(t)u_0 \in L_2(D) \right]$$

is $\left( \mathfrak{B}(\mathbb{R}^+) \otimes \mathfrak{F} \otimes \mathfrak{B}(L_2(D)), \mathfrak{B}(L_2(D)) \right)$-measurable.

In order to check that $\Phi^*$ is indeed the dual of $\Phi$, observe that

$$\begin{aligned}
\langle u_0, U_\omega^*(t)u_0^* \rangle &= \langle u_0, u^*(0, -t; \omega, u_0^*) \rangle \qquad &\text{(by (3.17))} \\
&= \langle u(t; 0, \theta_{-t}\omega, u_0), u_0^* \rangle \qquad &\text{(by Lemma 3.4)} \\
&= \langle U_{\theta_{-t}\omega}(t)u_0, u_0^* \rangle.
\end{aligned}$$

$\square$

We will call $\Phi$ as above the *measurable linear skew-product semiflow* on $L_2(D)$ *generated by* (3.1)+(3.2). The above construction of the measurable linear skew-product semiflow on $L_2(D)$, as well as its dual, can be repeated for the case when the zero-order term $c_0(\cdot, \cdot)$ is put to be constantly equal to zero, that is, for the problem

$$\partial_t u = \partial_i \big( a_{ij}(\theta_t \omega, x) \partial_j u + a_i(\theta_t \omega, x) u \big) + b_i(\theta_t \omega, x) \partial_i u, \quad t > s, \ x \in D, \qquad (3.18)$$

complemented with boundary condition

$$\mathcal{B}_\omega u = 0, \quad t > s, \ x \in \partial D, \qquad (3.19)$$

where $\mathcal{B}_\omega$ is the same as in (3.2), and its adjoint

$$-\partial_s u = \partial_i \big( a_{ji}(\theta_s \omega, x) \partial_j u - b_i(\theta_s \omega, x) u \big) - a_i(\theta_s \omega, x) \partial_i u, \quad s < t, \ x \in D, \qquad (3.20)$$

complemented with boundary condition

$$\mathcal{B}_\omega^* u = 0, \quad s < t, \ x \in \partial D, \qquad (3.21)$$

where $\mathcal{B}_\omega^*$ is the same as in (3.4).

For $\omega \in \Omega$, $s \in \mathbb{R}$ and $u_0 \in L_2(D)$ let $u^0(\cdot; s, \omega, u_0)$ stand for the global weak solution of $(3.18)_\omega + (3.19)_\omega$ satisfying the initial condition $u(s) = u_0$. Similarly, for $\omega \in \Omega$, $t \in \mathbb{R}$ and $u_0 \in L_2(D)$ let $u^{0*}(\cdot; t, \omega, u_0)$ stand for the global weak solution of $(3.20)_\omega + (3.21)_\omega$ satisfying the final condition $u(t) = u_0$. We write

$$U_\omega^0(t)u_0 := u^0(t, 0; \omega, u_0), \qquad (\omega \in \Omega, \ t \ge 0, \ u_0 \in L_2(D)), \qquad (3.22)$$

$$U_\omega^{0*}(t)u_0^* := u^{0*}(-t, 0; \omega, u_0^*), \qquad (\omega \in \Omega, \ t \ge 0, \ u_0^* \in L_2(D)). \qquad (3.23)$$

Since the coefficients of (3.18) + (3.19) and of (3.20) + (3.21) are bounded uniformly in $\omega$ and $x$, one obtains the following exponential estimate (see [4, Theorem 5.1]).

**Proposition 3.4** *There exists $\gamma \in \mathbb{R}$ such that*

$$\|U_\omega^0(t)\| \le e^{\gamma t} \tag{3.24}$$

*for all $\omega \in \Omega$ and $t > 0$.*

**Proposition 3.5** *For any $\omega \in \Omega$, $u_0 \in L_2(D)^+$ and $t > 0$,*

$$\exp\left(\int_0^t c_0^{(-)}(\theta_\tau \omega)\, d\tau\right) U_\omega^0(t) u_0 \le U_\omega(t) u_0 \le \exp\left(\int_0^t c_0^{(+)}(\theta_\tau \omega)\, d\tau\right) U_\omega^0(t) u_0.$$

*Proof* It follows by arguments as in the proof of [17, Lemma 4.3.1] that the solution of $(3.1)_\omega$ + $(3.2)_\omega$ with $c_0(\theta_t \omega, x)$ replaced by $c_0^{(\pm)}(\theta_t \omega)$, satisfying the initial condition $u(0, \cdot) = u_0$, equals

$$\left[ t \mapsto \exp\left(\int_0^t c_0^{(\pm)}(\theta_\tau \omega)\, d\tau\right) U_\omega^0(t) u_0 \right].$$

It suffices now to apply Proposition 3.2. □

### 3.2 Generalized Floquet Subspaces, Lyapunov Exponent, and Exponential Separation

In this subsection, we investigate the existence of generalized Floquet subspaces, Lyapunov exponent, and exponential separation. Throughout this subsection, we assume (PA0)–(PA3). We first introduce some further assumptions on $D$ and the coefficients of the problem (3.1) + (3.2).

**(PA5)** (Focusing) *There exist $\mathbf{e}, \mathbf{e}^* \in L_2(D)^+$, $\|\mathbf{e}\| = \|\mathbf{e}^*\| = 1$, $\tilde{\varkappa}, \tilde{\varkappa}^* \ge 1$, and $\tilde{\nu}, \tilde{\nu}^* > 0$ with the property that for each $\omega \in \Omega$ and any nonzero $u_0, u_0^* \in L_2(D)^+$ one can find $\tilde{\beta}(\omega, u_0), \tilde{\beta}^*(\omega, u_0^*) > 0$ such that*

$$\tilde{\beta}(\omega, u_0)\mathbf{e} \le U_\omega^0(1) u_0 \le \tilde{\varkappa}\tilde{\beta}(\omega, u_0)\mathbf{e}, \tag{3.25}$$
$$\tilde{\beta}^*(\omega, u_0^*)\mathbf{e}^* \le U_\omega^{0*}(1) u_0^* \le \tilde{\varkappa}^*\tilde{\beta}^*(\omega, u_0^*)\mathbf{e}^*. \tag{3.26}$$

*Moreover,*

$$\tilde{\beta}(\omega, \mathbf{e}) \ge \tilde{\nu}, \tag{3.27}$$
$$\tilde{\beta}^*(\omega, \mathbf{e}^*) \ge \tilde{\nu}^* \tag{3.28}$$

*for all $\omega \in \Omega$.*

At the end of this section, we give two sets of sufficient assumptions on the first and second-order coefficients, (R)(i), and (R)(ii), for the satisfaction of (PA4) and (PA5).
**(PA6)** (Zero order terms)

(i) *The mapping $\left[ \Omega \ni \omega \mapsto \int_0^1 c_0^{(+)}(\theta_t \omega)\, dt \in \mathbb{R}^+ \right]$ belongs to $L_1((\Omega, \mathfrak{F}, \mathbb{P}))$;*
(ii) *the mapping $\left[ \Omega \ni \omega \mapsto \ln^+ \int_0^1 \left(c_0^{(+)}(\theta_t \omega) - c_0^{(-)}(\theta_t \omega)\right) dt \in \mathbb{R}^+ \right]$ belongs to $L_1((\Omega, \mathfrak{F}, \mathbb{P}))$;*
(iii) *the mapping $\left[ \Omega \ni \omega \mapsto \int_0^1 \left(c_0^{(+)}(\theta_t \omega) - c_0^{(-)}(\theta_t \omega)\right) dt \in \mathbb{R}^+ \right]$ belongs to $L_1((\Omega, \mathfrak{F}, \mathbb{P}))$;*
(iv) *the mapping $\left[ \Omega \ni \omega \mapsto \int_0^1 c_0^{(-)}(\theta_t \omega)\, dt \in \mathbb{R}^- \right]$ belongs to $L_1((\Omega, \mathfrak{F}, \mathbb{P}))$.*

In the rest of this subsection, $\Phi = ((U_\omega(t)), (\theta_t))$ denotes the measurable linear skew-product semiflow on $L_2(D)$ generated by (3.1) + (3.2), and $\Phi^* = ((U_\omega^*(t)), (\theta_{-t}))$ denotes the dual of $\Phi$. The following are the main theorems of this subsection. Recall that we assume that (PA0) through (PA3) are fulfilled.

**Theorem 3.1** (Entire positive solution) *Assume* (PA4), (PA6)(i), *and that* (C1)(ii) *holds for* $\Phi$. *Moreover, assume that*

$$\limsup_{t \to \infty} \frac{\ln \|U_\omega(t)\|}{t} > -\infty.$$

*Then for* $\mathbb{P}$-*a.e.* $\omega \in \Omega$ *there exists a nontrivial entire positive solution of* (3.1)$_\omega$ + (3.2)$_\omega$.

**Theorem 3.2** (Generalized principal Floquet subspaces and Lyapunov exponent) *Assume* (PA4)–(PA5), *and* (PA6)(i)–(ii). *Then there are:*

- *an invariant set* $\tilde{\Omega}_0 \subset \Omega$, $\mathbb{P}(\tilde{\Omega}_0) = 1$,
- *an* $(\mathfrak{F}, \mathfrak{B}(L^2(D)))$-*measurable function* $w \colon \tilde{\Omega}_0 \to L_2(D)^+$ *with* $\|w(\omega)\| = 1$ *for all* $\omega \in \tilde{\Omega}_0$,
- *an* $(\mathfrak{F}, \mathfrak{B}(L^2(D)))$-*measurable function* $w^* \colon \tilde{\Omega}_0 \to L_2(D)^+$ *with* $\|w^*(\omega)\| = 1$ *for all* $\omega \in \tilde{\Omega}_0$,

*having the following properties:*

(i)
$$w(\theta_t \omega) = \frac{U_\omega(t) w(\omega)}{\|U_\omega(t) w(\omega)\|}$$

*for any* $\omega \in \tilde{\Omega}_0$ *and* $t \geq 0$.

(i)*
$$w^*(\theta_{-t} \omega) = \frac{U_\omega^*(t) w^*(\omega)}{\|U_\omega^*(t) w^*(\omega)\|}$$

*for any* $\omega \in \tilde{\Omega}_0$ *and* $t \geq 0$.

(ii) *Let for some* $\omega \in \tilde{\Omega}_0$ *a function* $v_\omega \colon \mathbb{R} \to L_2(D)^+ \setminus \{0\}$ *be an entire orbit of* $U_\omega$. *Then* $v_\omega(t) = \|v_\omega(0)\| w_\omega(t)$ *for all* $t \in \mathbb{R}$, *where*

$$w_\omega(t) := \begin{cases} \left(U_{\theta_t \omega}(-t)|_{\tilde{E}_1(\theta_t \omega)}\right)^{-1} w(\omega) & \text{for } t < 0 \\ U_\omega(t) w(\omega) & \text{for } t \geq 0, \end{cases}$$

*with* $\tilde{E}_1(\omega) = \mathrm{span}\{w(\omega)\}$.

(ii)* *Let for some* $\omega \in \tilde{\Omega}_0$ *a function* $v_\omega^* \colon \mathbb{R} \to L_2(D)^+ \setminus \{0\}$ *be an entire orbit of* $U_\omega^*$. *Then* $v_\omega^*(t) = \|v_\omega^*(0)\| w_\omega^*(t)$ *for all* $t \in \mathbb{R}$, *where*

$$w_\omega^*(t) := \begin{cases} \left(U_{\theta_{-t} \omega}^*(-t)|_{\tilde{E}_1^*(\theta_{-t} \omega)}\right)^{-1} w^*(\omega) & \text{for } t < 0 \\ U_\omega^*(t) w^*(\omega) & \text{for } t \geq 0, \end{cases}$$

*with* $\tilde{E}_1^*(\omega) = \mathrm{span}\{w^*(\omega)\}$.

(iii) *There is* $\tilde{\lambda}_1 \in [-\infty, \infty)$ *such that for any* $\omega \in \tilde{\Omega}_0$,

$$\tilde{\lambda}_1 = \lim_{t \to \pm\infty} \frac{\rho_t(\omega)}{t} = \int_\Omega \ln \rho_1 \, d\mathbb{P}$$

$$= \lim_{t \to \pm\infty} \frac{\ln \rho_t^*(\omega)}{t} = \int_\Omega \ln \rho_1^* \, d\mathbb{P},$$

*where*

$$\rho_t(\omega) = \begin{cases} \|U_\omega(t)w(\omega)\| & \text{for } t \geq 0 \\ 1/\|U_{\theta_t\omega}(-t)w(\theta_t\omega)\| & \text{for } t < 0, \end{cases}$$

$$\rho_t^*(\omega) = \begin{cases} \|U_\omega^*(t)w(\omega)\| & \text{for } t \geq 0 \\ 1/\|U_{\theta_{-t}\omega}^*(-t)w^*(\theta_{-t}\omega)\| & \text{for } t < 0. \end{cases}$$

(iv) *For any* $\omega \in \tilde{\Omega}_0$ *and* $u \in L_2(D)^+\backslash\{0\}$,

$$\limsup_{t\to\infty} \frac{\ln \|U_\omega(t)u\|}{t} = \tilde{\lambda}_1.$$

(v) *For any* $\omega \in \tilde{\Omega}_0$ *and* $u \in L_2(D)\backslash\{0\}$

$$\limsup_{t\to\infty} \frac{\ln \|U_\omega(t)u\|}{t} \leq \tilde{\lambda}_1$$

*and*

$$\lim_{t\to\infty} \frac{\ln \|U_\omega(t)\|}{t} = \tilde{\lambda}_1.$$

(vi) *Assume, moreover,* (PA6)(iv) *hold. Then* $\tilde{\lambda}_1 > -\infty$.

**Theorem 3.3** (Generalized exponential separation) *Assume* (PA4)–(PA5) *and* (PA6)(i),(iii). *Then the family* $\{\tilde{P}(\omega)\}_{\omega\in\tilde{\Omega}_0}$ *of projections associated with the invariant decomposition* $\tilde{E}_1(\omega) \oplus \tilde{F}_1(\omega) = L_2(D)$, *where* $\tilde{F}_1(\omega) = \{u_0 \in L_2(D) : \langle u_0, w^*(\omega)\rangle = 0\}$, *is strongly measurable and tempered. Moreover, the following holds.*

(i) $\tilde{F}_1(\omega) \cap L_2(D)^+ = \{0\}$ *for any* $\omega \in \tilde{\Omega}_0$.

(ii)

$$\lim_{t\to\infty} \frac{\ln \|U_\omega(t)\|}{t} = \lim_{t\to\infty} \frac{\ln \|U_\omega(t)u_0\|}{t} = \tilde{\lambda}_1$$

*for any* $\omega \in \tilde{\Omega}_0$ *and any* $u_0 \in L_2(D)\backslash\tilde{F}_1(\omega)$ *(in particular, for any nonzero* $u_0 \in L_2(D)^+$).

(iii) *There exist* $\tilde{\sigma} \in (0, \infty]$ *and* $\tilde{\lambda}_2 \in [-\infty, \infty)$, $\tilde{\lambda}_2 = \tilde{\lambda}_1 - \tilde{\sigma}$, *such that*

$$\lim_{t\to\infty} \frac{1}{t} \ln \frac{\|U_\omega(t)|_{\tilde{F}(\omega)}\|}{\|U_\omega(t)w(\omega)\|} = -\tilde{\sigma}$$

*and*

$$\lim_{t\to\infty} \frac{\ln \|U_\omega(t)|_{\tilde{F}_1(\omega)}\|}{t} = \tilde{\lambda}_2$$

*for each* $\omega \in \tilde{\Omega}_0$.

Before we prove Theorems 3.1–3.3, we first prove some propositions.

**Proposition 3.6** (Integrability) *Assume* (PA6)(i). *Then* (C1)(i) *and* (C1)$^*$(i) *hold for* $\Phi$ *and* $\Phi^*$, *respectively.*

*Proof* We prove the corresponding properties for $\Phi$ only, proofs for $\Phi^*$ being similar.

Since the norm on $L_2(D)$ is monotonic, it follows by Proposition 3.5 that

$$\|U_\omega(t)u_0\| \leq \exp\left(\int_0^t c_0^{(+)}(\theta_\tau\omega)\,d\tau\right)\|U_\omega^0(t)u_0\|$$

for any $\omega \in \Omega, t > 0$ and $u_0 \in L_2(D)^+$. As, by the Banach lattice property, each $u_0 \in L_2(D)$ can be written as $u_0^+ - u_0^-$, with $u_0^\pm \in L_2(D)^+$, $\|u_0^\pm\| \leq \|u_0\|$, we have

$$\|U_\omega(t)\| \leq 2 \exp \left( \int_0^t c_0^{(+)}(\theta_\tau \omega) \, d\tau \right) \|U_\omega^0(t)\|.$$

Proposition 3.4 together with (PA6)(i) conclude the proof. $\qquad \square$

*Remark 3.1* (1) (Compactness) By Proposition 3.1 (iv), (C1)(iii) and (C1)*(iii) are satisfied for $\Phi$ and $\Phi^*$.
(2) (Positivity) By Proposition 3.1 (iii), (C2) and (C2)* are satisfied for $\Phi$ and $\Phi^*$.

**Proposition 3.7** (Focusing) *Assume* (PA5) *and* (PA6)(ii)*. Then* (C3) *and* (C3)* *are satisfied for* $\Phi$ *and* $\Phi^*$*, respectively.*

*Proof* We prove (C3). (C3)* can be proved similarly.

Proposition 3.5 together with (PA5) implies that

$$\beta(u_0, \omega)\mathbf{e} \leq U_\omega(1)u_0 \leq \varkappa(\omega)\beta(u_0, \omega)\mathbf{e} \quad \forall \, u_0 \in L_2(D)^+,$$

where $\beta(u_0, \omega) = \tilde{\beta}(\omega, u_0) \exp(\int_0^1 c_0^{(-)}(\theta_t \omega) \, dt)$ and $\varkappa(\omega) = \tilde{\varkappa} \exp(\int_0^1 (c_0^{(+)}(\theta_t \omega) - c_0^{(-)}(\theta_t \omega)) \, dt)$. This together with (PA6)(ii) implies (C3). $\qquad \square$

**Proposition 3.8** (Strong focusing) *Assume* (PA5) *and* (PA6)(iii)*. Then* (C4) *holds.*

*Proof* First of all, by Proposition 3.7, (C3) and (C3)* are satisfied for $\Phi$ and $\Phi^*$, respectively.

Next, by the arguments of Proposition 3.7,

$$\varkappa(\omega) = \tilde{\varkappa} \exp \left( \int_0^1 \left( c_0^{(+)}(\theta_t \omega) - c_0^{(-)}(\theta_t \omega) \right) dt \right)$$

and

$$\varkappa^*(\omega) = \tilde{\varkappa}^* \exp \left( \int_0^1 \left( c_0^{(+)}(\theta_t \omega) - c_0^{(-)}(\theta_t \omega) \right) dt \right).$$

Then by (PA6)(iii),

$$\ln \varkappa, \ \ln \varkappa^* \in L_1((\Omega, \mathfrak{F}, \mathbb{P})).$$

Now, by [17, Proposition 2.2.9(2)], for any $u_0, u_0^* \in L_2(D)^+ \backslash \{0\}$,

$$\left( U_\omega^0(1)u_0 \right)(x) > 0 \quad \forall \, x \in D$$

and

$$\left( U_\omega^{0*}(1)u_0^* \right)(x) > 0 \quad \forall \, x \in D.$$

Then by (PA5),

$$\mathbf{e}(x) > 0 \quad \text{and} \quad \mathbf{e}^*(x) > 0 \quad \forall \, x \in D$$

and hence

$$\langle \mathbf{e}, \mathbf{e}^* \rangle > 0.$$

Therefore (C4) holds. $\qquad \square$

**Proposition 3.9** (Strong positivity in one direction) *Assume* (PA5) *and* (PA6)(iv). *Then* (C5) *and* (C5)* *hold.*

*Proof* We prove the corresponding properties for $\Phi$ only, proofs for $\Phi^*$ being similar.

By (3.25), (3.27) and Proposition 3.5,

$$U_\omega(1)\mathbf{e} \geq \tilde{\nu} \exp\left(\int_0^1 c_0^{(-)}(\theta_t\omega)\,dt\right)\mathbf{e}$$

for any $\omega \in \Omega$. This together with (PA6)(iv) implies that (C5) holds. □

We now prove Theorems 3.1–3.3. Observe that (C0) and (C0)* are satisfied for $X = X^* = L_2(D)$.

*Proof of Theorem 3.1* Assume (PA4), (PA6)(i), and that (C1)(ii) hold for $\Phi$. By Proposition 3.6 and Remark 3.1, (C1)(i)–(iii) and (C2) are satisfied. By Proposition 2.3,

$$\lambda_1 = \lim_{t\to\infty} \frac{\ln\|U_\omega(t)\|}{t} > -\infty \quad \text{for } \mathbb{P}\text{-a.e. } \omega \in \Omega.$$

The theorem then follows from Theorem 2.1. □

*Proof of Theorem 3.2* Assume (PA4)–(PA5), and (PA6)(i)–(ii). By Proposition 3.6 and Remark 3.1, (C1)(i), (C1)*(i), (C2), and (C2)* are satisfied. By Proposition 3.7, (C3) and (C3)* are satisfied. Theorem 3.2(i)–(v) then follows from Theorem 2.2(1)–(5) and Theorem 2.3.

Assume, moreover, (PA6)(iv) holds. Then by Proposition 3.9, (C5) holds. Theorem 3.2(vi) then follows from Theorem 2.2(7). □

*Proof of Theorem 3.3* Assume (PA4)–(PA5) and (PA6)(i),(iii). By Proposition 3.6 and Remark 3.1, (C1)(i), (C1)*, (C2), and (C2)* are satisfied. By Proposition 3.8, (C4) is satisfied. Theorem 3.3 then follows from Theorem 2.4. □

We now give sufficient conditions, **(R)(i)**, **(R)(ii)**, for the satisfaction of (PA4) and (PA5).

**(R)** *Either of the assumptions below*, (R)(i) or (R)(ii), *is satisfied:*
**(R)(i)** (Higher-order coefficients independent of $\omega$) *In the Dirichlet boundary condition case:*

- *The functions $a_{ij}$ $(i, j = 1, \ldots, N)$, $a_i$ $(i = 1, \ldots, N)$, $b_i$ $(i = 1, \ldots, N)$ depend on $x$ only, and belong to $L_\infty(D)$.*

In the Robin boundary condition case:

- *$D$ is a bounded domain, where its boundary $\partial D$ is an $(N-1)$-dimensional manifold of class $C^2$.*
- *The functions $a_{ij}$ $(i, j = 1, \ldots, N)$, $a_i$ $(i = 1, \ldots, N)$, $b_i$ $(i = 1, \ldots, N)$ depend on $x$ only, and belong to $C^1(\bar{D})$.*
- *The function $d_0$ depends on $x$ only, and belongs to $C^1(\partial D)$.*

**(R)(ii)** *(Classical case) There exists $\alpha \in (0, 1)$ such that*

- *$\partial D$ is an $(N-1)$-dimensional manifold of class $C^{3+\alpha}$.*

- *For each $\omega \in \Omega$ the functions $a_{ij}^{\omega}$ ($i, j = 1, \ldots, N$) and $a_i^{\omega}$ ($i = 1, \ldots, N$) belong to $C^{2+\alpha,2+\alpha}(\mathbb{R} \times \bar{D})$. Moreover, their $C^{2+\alpha,2+\alpha}(\mathbb{R} \times \bar{D})$-norms are bounded uniformly in $\omega \in \Omega$.*
- *For each $\omega \in \Omega$ the functions $b_i^{\omega}$ ($i = 1, \ldots, N$) belong to $C^{2+\alpha,1+\alpha}(\mathbb{R} \times \bar{D})$. Moreover, their $C^{2+\alpha,1+\alpha}(\mathbb{R} \times \bar{D})$-norms are bounded uniformly in $\omega \in \Omega$.*
- *In the Robin boundary condition case, for each $\omega \in \Omega$ the function $d_0^{\omega}$ belongs to $C^{2+\alpha,2+\alpha}(\mathbb{R} \times \partial D)$. Moreover, their $C^{2+\alpha,2+\alpha}(\mathbb{R} \times \partial D)$-norms are bounded uniformly in $\omega \in \Omega$.*

**Proposition 3.10** *If* (R) *holds, then* (PA4) *and* (PA5) *are satisfied.*

*Proof* First, assume (R)(i). Regarding (PA4) there is nothing to prove. As the coefficients are independent of $\omega$, we write the global solution of (3.1) + (3.2) (in the case $c_0 \equiv 0$) with initial value $u(0) = u_0 \in L_2(D)$ as $u^0(\cdot; u_0)$. By Proposition 3.1(iii),

$$u^0(1; u_0) \not\equiv 0 \quad \forall u_0 \in L_2(D)^+ \setminus \{0\}. \tag{3.29}$$

We claim that there are $\mathbf{e} \in L_2(D)^+$, $C > 0$, and $\tilde{\beta}: L_2(D)^+ \setminus \{0\} \to (0, \infty)$ such that

$$\tilde{\beta}(u_0)\mathbf{e} \le u^0(1; u_0) \le C\tilde{\beta}(u_0)\mathbf{e} \tag{3.30}$$

for any nonzero $u_0 \in L_2(D)^+$.

In fact, it follows from [8,9] that there is a simple eigenvalue $\lambda_{\text{princ}}$ (the *principal eigenvalue*) of the problem

$$\begin{cases} \partial_i \big(a_{ij}(x)\partial_j v + a_i(x)v\big) + b_i(x)\partial_i v = \lambda v, & x \in D, \\ \mathcal{B}_\omega v = 0, & x \in \partial D, \end{cases} \tag{3.31}$$

which is real and larger than and bounded away from the real parts of the remaining eigenvalues, and that an eigenfunction corresponding to it (a *principal eigenfunction*) can be chosen so to take positive values on $D$ (note that $\mathcal{B}_\omega$ in (3.31) is independent of $\omega$).

As $\mathbf{e}$ we take the principal eigenfunction, normalized so that $\|\mathbf{e}\| = 1$. In the Dirichlet boundary condition case, by [9, Theorem 2.1], there is a constant $C > 1$ such that

$$\sup_{x \in D} \frac{u^0(1; u_0)(x)}{u^0(1; \mathbf{e})(x)} \le C \inf_{x \in D} \frac{u^0(1; u_0)(x)}{u^0(1; \mathbf{e})(x)} \tag{3.32}$$

for any $u_0 \in L_2(D)^+ \setminus \{0\}$. In the Robin boundary condition case, by [8, Theorem 2.5], there is a constant $\bar{C} > 1$ such that

$$\sup_{x \in \bar{D}} u^0(1; u_0)(x) \le \bar{C} \inf_{x \in \bar{D}} u^0(1; u_0)(x) \quad \forall u_0 \in L_2(D)^+. \tag{3.33}$$

Then by [17, Lemma 3.3.1], (3.32) also holds.

As $u^0(1; \mathbf{e}) = e^{\lambda_{\text{princ}}} \mathbf{e}$, we obtain

$$\tilde{\beta}(u_0)\mathbf{e} \le u^0(1; u_0) \le C\tilde{\beta}(u_0)\mathbf{e}$$

for any nonzero $u_0 \in L_2(D)^+$, where

$$\tilde{\beta}(u_0) = e^{\lambda}_{\text{princ}} \inf_{x \in D} \frac{u^0(1; u_0)(x)}{u^0(1; \mathbf{e})(x)}$$

is positive, since otherwise, by (3.32), $u^0(1; u_0)$ would be constantly equal to zero, which contradicts (3.29). Therefore, (3.30) holds, which implies that (3.25) is satisfied. (3.26) is proved in an analogous way. Regarding (3.27) and (3.28), there is nothing to prove.

Second, we assume (R)(ii). We embed problem into a family of problems as in [17]. For $\omega \in \Omega$ put

$$a^{0,\omega} := \big((a_{ij}^{\omega})_{i,j=1}^N, (a_i^{\omega})_{i=1}^N, (b_i^{\omega})_{i=1}^N, 0, d_0^{\omega}\big)$$

(in the Dirichlet case we put $d_0^{\omega}$ constantly equal to zero). Let $Y^0$ be the closure of the set $\{a^{0,\omega} : \omega \in \Omega\}$ in the weak-* topology of $L_{\infty}(\mathbb{R} \times D, \mathbb{R}^{N^2+2N+1}) \times L_{\infty}(\mathbb{R} \times \partial D, \mathbb{R})$. We define $E^0 : \Omega \to Y^0$ by $E^0(\omega) := a^{0,\omega}$.

For each $\tilde{a}^0 = \big((\tilde{a}_{ij})_{i,j=1}^N, (\tilde{a}_i)_{i=1}^N, (\tilde{b}_i)_{i=1}^N, 0, \tilde{d}_0\big) \in Y^0$ we have that $\tilde{a}_{ij}, \tilde{a}_i \in C^{2+\alpha,2+\alpha}(\mathbb{R} \times \bar{D})$, $\tilde{b}_i \in C^{2+\alpha,1+\alpha}(\mathbb{R} \times \bar{D})$ and $\tilde{d}_0 \in C^{2+\alpha,2+\alpha}(\mathbb{R} \times \bar{D})$, with the corresponding norms bounded uniformly in $\tilde{a}^0 \in Y^0$. Also, on $Y^0$ the weak-* and open-compact topologies coincide. Consequently, (PA4) holds.

We consider, for each $s \in \mathbb{R}$ and each $\tilde{a}^0 \in Y^0$,

$$\partial_t u = \partial_i\big(\tilde{a}_{ij}(t,x)\partial_j u + \tilde{a}_i(t,x)u\big) + \tilde{b}_i(t,x)\partial_i u, \quad t \in (s,\infty), \ x \in D \qquad (3.34)$$

complemented with boundary conditions

$$\mathcal{B}_{\tilde{a}^0} u = 0, \quad t \in (s,\infty), \ x \in \partial D, \qquad (3.35)$$

where $\mathcal{B}_{\tilde{a}^0} = \mathcal{B}$ in the Dirichlet or periodic cases, or

$$\mathcal{B}_{\tilde{a}^0} u = \big(\tilde{a}_{ij}(t,x)\partial_j u + \tilde{a}_i(t,x)u\big)v_i + \tilde{d}_0(t,x)u \qquad (3.36)$$

in the Robin case.

For $s \in \mathbb{R}$, $\tilde{a}^0 \in Y^0$ and $u_0 \in L_2(D)$ let $\tilde{u}^0(\cdot; s, \tilde{a}^0, u_0)$ stand for the global weak solution of $(3.34)_{\tilde{a}^0} + (3.35)_{\tilde{a}^0}$ satisfying the initial condition $u(s) = u_0$. By the linearity of the problems, we can (and do) write $U_{\tilde{a}^0}^0(t)u_0$ for $\tilde{u}^0(t; 0, \tilde{a}^0, u_0)$.

By [17, Theorem 3.3.1], there is a continuous function $w^0 : Y^0 \to L_2(D)^+$ with $\|w^0(\tilde{a}^0)\| = 1$ for any $\tilde{a}^0 \in Y^0$ such that

$$U_{\tilde{a}^0}^0(t)w^0(\tilde{a}^0) = \|U_{\tilde{a}^0}^0(t)w^0(\tilde{a}^0)\| w^0(\tilde{a}^0 \cdot t) \quad \forall t \geq 0, \qquad (3.37)$$

where $\tilde{a}^0 \cdot t$ denotes the $t$-translate of $\tilde{a}^0$. In the Dirichlet boundary condition case, let $\mathbf{e}$ be the positive principal eigenfunction of

$$\begin{cases} \Delta u = \lambda u, & x \in D \\ u = 0, & x \in \partial D \end{cases}$$

with $\|\mathbf{e}\| = 1$ and in the Robin boundary condition case, let $\mathbf{e}$ be a constant positive function with $\|\mathbf{e}\| = 1$. We claim that there are positive constants $C_-, C_+ > 0$ such that

$$C_- \mathbf{e} \leq w^0(\tilde{a}^0) \leq C_+ \mathbf{e} \quad \forall \tilde{a}^0 \in Y^0. \qquad (3.38)$$

In fact, assume that $C_-$ does not exist. Then there is a sequence $(\tilde{a}^{(n)}) \subset Y^0$ such that

$$w^0(\tilde{a}^{(n)}) \geq \frac{1}{n}\mathbf{e}$$

does not hold for all $n \geq 1$. By (R)(ii), without loss of generality, we may assume that $\tilde{a}^{(n)} \to \tilde{a}^0 \in Y^0$ as $n \to \infty$. This implies that

$$w^0(\tilde{a}^{(n)}) \to w^0(\tilde{a}^0)$$

as $n \to \infty$ in $L^2(D)$. [17, Lemma 3.5.1] implies that

$$w^0(\tilde{a}^{(n)}) \to w^0(\tilde{a}^0)$$

in $C^1(\bar{D})$. As, by (3.37), $w^0(\tilde{a}^0)$ equals $U^0_{\tilde{a}^0\cdot(-1)}(1)w^0(\tilde{a}^0 \cdot (-1))$ multiplied by a positive number and $w^0(\tilde{a}^0\cdot(-1)) \in L_2(D)^+\backslash\{0\}$, from the Hopf maximum principle (in the Dirichlet case) or the strong maximum principle (in the Neumann or Robin cases) we deduce that there is $C_0 > 0$ such that $w^0(\tilde{a}^0) > C_0\mathbf{e}$, so

$$w^0(\tilde{a}^{(n)}) \geq C_0\mathbf{e} \quad \forall n \gg 1.$$

This is a contradiction. Hence $C_-$ exists. Similarly, we can prove that $C_+$ exists.

We also claim that there is $C > 0$ such that

$$\sup_{x\in D}\frac{\left(U^0_{\tilde{a}^0}(1,0)u_0\right)(x)}{\left(U^0_{\tilde{a}^0}(1,0)w(\tilde{a}^0)\right)(x)} \leq C \inf_{x\in D}\frac{\left(U^0_{\tilde{a}^0}(1,0)u_0\right)(x)}{\left(U^0_{\tilde{a}^0}(1,0)w(\tilde{a}^0)\right)(x)} \tag{3.39}$$

for all $\tilde{a}^0 \in Y^0$ and $u_0 \in L_2(D)^+$. In fact, in the Dirichlet boundary condition case, (3.39) follows from [9, Theorem 2.1]. In the Neumann or Robin boundary condition cases, by [8, Theorem 2.5], there is a constant $\tilde{C} > 1$ such that

$$\sup_{x\in\bar{D}}\left(U_{\tilde{a}^0}(1,0)u_0\right)(x) \leq \tilde{C} \inf_{x\in\bar{D}}\left(U_{\tilde{a}^0}(1,0)w(\tilde{a}^0)\right)(x) \quad \forall u_0 \in L_2(D)^+, \ \tilde{a}^0 \in Y^0.$$

Then by [17, Lemma 3.3.1], (3.39) also holds.

By (3.38) and (3.39),

$$\tilde{\beta}(u_0)\mathbf{e} \leq U^0_{\tilde{a}^0}(1,0)u_0 \leq C\tilde{\beta}(u_0)\mathbf{e} \quad \forall u_0 \in L_2(D)^+, \tag{3.40}$$

where

$$\tilde{\beta}(u_0) = \inf_{x\in D}\frac{\left(U_{\tilde{a}^0}(1,0)u_0\right)(x)}{\left(U_{\tilde{a}^0}(1,0)w(\tilde{a}^0)\right)(x)}$$

is positive, since otherwise, by (3.39), $U_{\tilde{a}^0}(1,0)u_0$ would be constantly equal to zero. This implies (3.25).

To prove (3.27), suppose to the contrary that there is a sequence $(\tilde{a}^{(n)}) \subset Y^0$ such that

$$U_{\tilde{a}^{(n)}}(1,0)\mathbf{e} \geq \frac{1}{n}\mathbf{e}$$

does not hold for all $n \geq 1$. By (R)(ii), without loss of generality, we may assume that $\tilde{a}^{(n)} \to \tilde{a}^0 \in Y^0$ as $n \to \infty$. [17, Proposition 2.5.4] implies that

$$U_{\tilde{a}^{(n)}}(1,0)\mathbf{e} \to U_{\tilde{a}^0}(1,0)\mathbf{e}$$

in $C^1(\bar{D})$. By the Hopf maximum principle (in the Dirichlet case) or the strong maximum principle (in the Neumann or Robin cases), we can find $\hat{C} > 0$ such that $U_{\tilde{a}^0}(1,0)\mathbf{e} > \hat{C}\mathbf{e}$, so

$$U_{\tilde{a}^{(n)}}(1,0)\mathbf{e} \geq C_0\mathbf{e} \quad \forall n \gg 1.$$

This is a contradiction. Hence (3.27) is satisfied.

The fulfillment of (3.26) and (3.28) follows by applying analogous reasoning to the adjoint problem.                                                                                              □

## 4 Systems of Linear Random Delay Differential Equations

In this section, we consider applications of the general results stated in Section 2 to (1.3), that is, the following system of linear random delay differential equations,

$$\frac{du}{dt} = A(\theta_t\omega)u(t) + B(\theta_t\omega)u(t-1) \tag{4.1}$$

where $u \in \mathbb{R}^N$, $N \geq 2$, and $A(\omega)$, $B(\omega)$ are $N$ by $N$ real matrices (we write $A(\omega)$, $B(\omega) \in \mathbb{R}^{N \times N}$):

$$A(\omega) = \begin{pmatrix} a_{11}(\omega) & a_{12}(\omega) & \cdots & a_{1N}(\omega) \\ a_{21}(\omega) & a_{22}(\omega) & \cdots & a_{2N}(\omega) \\ \vdots & \vdots & \ddots & \vdots \\ a_{N1}(\omega) & a_{N2}(\omega) & \cdots & a_{NN}(\omega) \end{pmatrix},$$

and

$$B(\omega) = \begin{pmatrix} b_{11}(\omega) & b_{12}(\omega) & \cdots & b_{1N}(\omega) \\ b_{21}(\omega) & b_{22}(\omega) & \cdots & b_{2N}(\omega) \\ \vdots & \vdots & \ddots & \vdots \\ b_{N1}(\omega) & b_{N2}(\omega) & \cdots & b_{NN}(\omega) \end{pmatrix}.$$

Again, let $((\Omega, \mathfrak{F}, \mathbb{P}), (\theta_t)_{t \in \mathbb{R}})$ be an ergodic metric dynamical system, with $\mathbb{P}$ complete.

In this section, as a Banach space $X$ we will consider the space $C([-1, 0], \mathbb{R}^N)$ of continuous $\mathbb{R}^N$-valued functions defined on $[-1, 0]$.

The symbol $\|\cdot\|$ stands, depending on the context, either for the Euclidean norm on $\mathbb{R}^N$, or for the Euclidean matrix norm on the algebra $\mathbb{R}^{N \times N}$, or else for the corresponding maximum norm on $C([-1, 0], \mathbb{R}^N)$:

$$\|u\| = \max_{t \in [-1,0]} \|u(t)\| \quad \text{for} \quad u(\cdot) \in C([-1, 0], \mathbb{R}^N).$$

Similarly, $\|\cdot\|_1$ stands either for the $\ell_1$-norm on $\mathbb{R}^N$ or for the corresponding maximum norm on $C([-1, 0], \mathbb{R}^N)$:

$$\|u\|_1 = \max_{t \in [-1,0]} \|u(t)\|_1 \quad \text{for} \quad u(\cdot) \in C([-1, 0], \mathbb{R}^N).$$

We will use the notation $\leq$ (and $\geq$) to denote the order relations generated by the standard cone $(\mathbb{R}^N)^+$ in $\mathbb{R}^N$ as well as the standard cone $C([-1, 0], \mathbb{R}^N)^+$ in $C([-1, 0], \mathbb{R}^N)$.

Throughout this section, we make the following standing assumption.

**(OA0)** (Measurability) $A, B \colon \Omega \to \mathbb{R}^{N \times N}$ are $(\mathfrak{F}, \mathfrak{B}(\mathbb{R}^{N \times N}))$-measurable, and $[\, \mathbb{R} \ni t \mapsto A(\theta_t\omega) \in \mathbb{R}^{N \times N} \,]$, $[\, \mathbb{R} \ni t \mapsto B(\theta_t\omega) \in \mathbb{R}^{N \times N} \,]$ *are continuous for all* $\omega \in \Omega$.

Under the assumption (OA0), for any $u_0 \in C([-1, 0], \mathbb{R}^N)$ and any $\omega \in \Omega$, there is a unique function $\big[\, [-1, \infty) \ni t \mapsto u(t; \omega, u_0) \in \mathbb{R}^N \,\big]$ such that

- Equation (4.1) is satisfied for each $t \geq 0$, where for $t = 0$ we understand the right-hand derivative;
- the initial condition
$$u(\tau; \omega, u_0) = u_0(\tau) \quad \forall \, \tau \in [-1, 0] \tag{4.2}$$
holds.

For a proof see, e.g., [7, Chapter 2].

We give now a useful representation of the solution of (4.1)+(4.2). Namely, for $\omega \in \Omega$ and $u_0^0 \in \mathbb{R}^N$ denote by $U_\omega^0(t)u_0^0$ the value at time $t \geq 0$ of the solution of $u' = A(\theta_t\omega)u$ satisfying the initial condition $u(0) = u_0^0$. For each $\omega \in \Omega$ and each $u_0 \in C([-1, 0], \mathbb{R}^N)$ the function $u(\cdot; \omega, u_0)$ is the solution of the integral equation

$$u(t; \omega, u_0) = U_\omega^0(t)u_0(0) + \int_0^t U_{\theta_\tau \omega}^0(t - \tau) B(\theta_\tau \omega) u(\tau - 1; \omega, u_0) \, d\tau, \quad t > 0. \quad (4.3)$$

For $t \geq 0$, $\omega \in \Omega$ and $u_0 \in C([-1, 0], \mathbb{R}^N)$ define $U_\omega(t)u_0 \in C([-1, 0], \mathbb{R}^N)$ by the formula

$$(U_\omega(t)u_0)(\tau) := u(t + \tau; \omega, u_0), \quad \tau \in [-1, 0].$$

**Proposition 4.1** *Under* (OA0)*,* $\left((U_\omega(t))_{\omega \in \Omega, t \in [0,\infty)}, (\theta_t)_{t \in \mathbb{R}}\right)$ *is a measurable linear skew-product semiflow on* $C([-1, 0], \mathbb{R}^N)$ *covering* $(\theta_t)$.

*Proof* The satisfaction of (2.1) and (2.2) is a standard exercise.

The continuity, for any $\omega \in \Omega$ and any $u_0 \in C([-1, 0], \mathbb{R}^N)$, of the mapping $[\, [0, \infty) \ni t \mapsto U_\omega(t)u_0 \in C([-1, 0], \mathbb{R}^N)\,]$ is straightforward. The fact that for each $\omega \in \Omega$ and $t \geq 0$ the mapping $U_\omega(t)$ is a bounded linear operator from $C([-1, 0], \mathbb{R}^N)$ into $C([-1, 0], \mathbb{R}^N)$ is a consequence of the continuous dependence of solutions on initial values (see, e.g., [7, Chapter 2, in particular Section 2.6]).

Fix now $u_0 \in C([-1, 0], \mathbb{R}^N)$ and $t > 0$. In view of Proposition 2.1 it is enough to prove the measurability of the mapping $\left[\, \Omega \ni \omega \mapsto U_\omega(t)u_0 \in C([-1, 0], \mathbb{R}^N)\,\right]$. Further, it follows from (2.2) that it suffices to prove that for $t \in (0, 1]$.

By a variant of Pettis' Measurability Theorem (see, e.g., [6, Corollary 4 on pp. 42–43]) combined with other results on measurable Banach-space-valued functions [1, Theorem 4.38 on p. 145 and Lemma 11.37 on p. 424], the mapping $\left[\, \Omega \ni \omega \mapsto U_\omega(t)u_0 \in C([-1, 0], \mathbb{R}^N)\,\right]$ is $(\mathfrak{F}, \mathfrak{B}\left(C([-1, 0], \mathbb{R}^N)\right))$-measurable if and only if the mapping $\left[\, \Omega \ni \omega \mapsto (U_\omega(t)u_0)(\tau) \in \mathbb{R}^N \,\right]$ is $(\mathfrak{F}, \mathfrak{B}(\mathbb{R}^N))$-measurable for each $\tau \in [-1, 0]$. Consequently, the problem reduces to proving the $(\mathfrak{F}, \mathfrak{B}(\mathbb{R}^N))$-measurability of the mapping $\left[\, \Omega \ni \omega \mapsto u(t; \omega, u_0) \in \mathbb{R}^N \,\right]$ for each $t \in (0, 1]$.

Observe that for such $t$ (4.3) takes the form

$$u(t; \omega, u_0) = U_\omega^0(t)u_0(0) + \int_0^t U_{\theta_\tau \omega}^0(t - \tau) B(\theta_\tau \omega) u_0(\tau - 1) \, d\tau.$$

By repeated application of the variant of Pettis' Measurability Theorem mentioned above, together with (OA0) and the fact that $\left[\, \Omega \times [0, 1] \ni (\omega, s) \mapsto U_\omega(s) \in \mathbb{R}^{N \times N} \,\right]$ is $(\mathfrak{F} \otimes \mathfrak{B}([0, 1]), \mathfrak{B}(\mathbb{R}^{N \times N}))$-measurable (see, e.g., [2, Example 2.2.8]) together with (OA0) we obtain the desired result. $\qquad\qquad\square$

From now until the end of the present section we assume that (OA0) holds. We will call $\left((U_\omega(t))_{\omega \in \Omega, t \in [0,\infty)}, (\theta_t)_{t \in \mathbb{R}}\right)$ the *measurable linear skew-product semiflow on* $C([-1, 0], \mathbb{R}^N)$ *generated by* (4.1).

**Proposition 4.2** (Compactness) *For any* $k \in \mathbb{N}$, $U_\omega(k)$ *is, for each* $\omega \in \Omega$, *completely continuous.*

*Proof* Let $k \in \mathbb{N}$ be given. Fix $\omega \in \Omega$. Since for each $u_0 \in C([-1, 0], \mathbb{R}^N)$ the mapping $\left[\, [0, k] \ni t \mapsto \|U_\omega(t)u_0\| \in [0, \infty) \,\right]$ is continuous, the Uniform Boundedness Theorem implies that the set $\{\,\|U_\omega(t)\| : t \in [0, k]\,\}$ is bounded. Denote its supremum by $M$. Obviously, $M \geq 1$. For any $k - 1 \leq t_1 < t_2 \leq k$ and any $u_0 \in C([-1, 0], \mathbb{R}^N)$ we estimate

$$\|u(t_2; \omega, u_0) - u(t_1; \omega, u_0)\|$$

$$\leq \int_{t_1}^{t_2} \|A(\theta_s\omega)\| \|u(s; \omega, u_0)\|\, ds + \int_{t_1}^{t_2} \|B(\theta_s\omega)\| \|u(s-1; \omega, u_0)\|\, ds$$

$$\leq \left( \int_{t_1}^{t_2} \max\{\|A(\theta_s\omega)\|, \|B(\theta_s\omega)\|\}\, ds \right) (M+1)\|u_0\|,$$

which implies, via (OA0), that the set $\{\, U_\omega(k)u_0 : \|u_0\| \leq 1\,\}$ is equicontinuous. $\qquad\square$

To investigate the generalized principal Floquet spaces and principal Lyapunov exponent of $\big((U_\omega(t))_{\omega\in\Omega, t\in[0,\infty)}, (\theta_t)_{t\in\mathbb{R}}\big)$, we state more standing assumptions on $A(\omega)$ and $B(\omega)$. Let

$$\underline{a}_{ii}(\omega) := \min_{\substack{0 \leq s < t \leq 1 \\ k=0,1,\dots,N+1}} \int_{k+s}^{k+t} a_{ii}(\theta_\tau\omega)\, d\tau \tag{4.4}$$

and

$$\underline{a}(\omega) := \min\{e^{\underline{a}_{11}(\omega)}, e^{\underline{a}_{22}(\omega)}, \dots, e^{\underline{a}_{NN}(\omega)}\}. \tag{4.5}$$

**(OA1)** (Cooperativity)

(i) $a_{ij}(\omega) \geq 0$ *for all* $i \neq j$, $i, j = 1, 2, \dots, N$ *and* $\omega \in \Omega$.
(ii) $b_{ij}(\omega) \geq 0$ *for all* $i, j = 1, 2, \dots, N$ *and* $\omega \in \Omega$.

**(OA2)** (Integrability)

(i) *The function* $[\,\Omega \ni \omega \mapsto \max_{1 \leq i,j \leq N} a_{ij}(\omega)\,]$ *is in* $L_1((\Omega, \mathfrak{F}, \mathbb{P}))$.
(ii) *The function* $[\,\Omega \ni \omega \mapsto \max_{1 \leq i,j \leq N} b_{ij}(\omega)\,]$ *is in* $L_1((\Omega, \mathfrak{F}, \mathbb{P}))$.

**(OA3)** (Irreducibility) *There is an* $(\mathfrak{F}, \mathfrak{B}(\mathbb{R}))$-*measurable function* $\underline{\delta} \colon \Omega \to (0, \infty)$ *such that for each* $\omega \in \Omega$ *and* $i \in \{1, 2, \dots, N\}$ *there are* $j_0 = i, j_2, j_3, \dots, j_{N-1} \in \{1, 2, \dots, N\}$ *satisfying*

(i) $\{j_0, j_1, \dots, j_{N-1}\} = \{1, 2, \dots, N\}$ *and* $b_{j_{l+1}j_l}(\theta_t\omega) \geq \underline{\delta}(\omega)$ *for* $0 \leq t \leq N + 2$ *and* $l = 0, 1, \dots, N - 2$.
(ii) $\ln^+ \ln\big(\overline{\beta}/\underline{\beta}\big) \in L_1((\Omega, \mathfrak{F}, \mathbb{P}))$, *where*

$$\underline{\beta}(\omega) = \min\left\{ \frac{\underline{a}^{N+2}(\omega)\underline{\delta}^k(\omega)}{k!}, \frac{\underline{a}^{N+3}(\omega)\underline{\delta}^k(\omega)}{k!} : k = 1, 2, \dots, N \right\},$$

$$\overline{\beta}(\omega) = \max\{\bar{a}^{N+2}(\omega)(1 + N\bar{\delta}(\omega))^{N+1}, \bar{a}^{N+2}(\omega)N\bar{\delta}(\omega)(1 + N\bar{\delta}(\omega))^{N+1}\},$$

$$\bar{a}(\omega) = \exp\left( \int_0^{N+2} \sum_{l=1}^{N} \max_{1 \leq k \leq N} a_{lk}(\theta_\tau\omega)\, d\tau \right), \quad \bar{\delta}(\omega) = \max_{\substack{0 \leq t \leq N+2 \\ i,j=1,2,\dots,N}} b_{ij}(\theta_t\omega).$$

(iii) $\ln^- \underline{\beta} \in L_1((\Omega, \mathfrak{F}, \mathbb{P}))$, *where* $\underline{\beta}$ *is as in (ii).*

**(OA4)** (Positivity)

(i) *There is an $(\mathfrak{F}, \mathfrak{B}(\mathbb{R}))$-measurable function $\underline{\delta}\colon \Omega \to (0, \infty)$ such that for any $\omega \in \Omega$ and any $i \neq j$ there holds $\max\{b_{ij}(\theta_t\omega)\} \geq \underline{\delta}(\omega)$ for $t \in [0, 2]$.*

(ii) $\ln^+ \ln\left(\overline{\overline{\beta}}/\underline{\underline{\beta}}\right) \in L_1((\Omega, \mathfrak{F}, \mathbb{P}))$, where

$$\underline{\underline{\beta}}(\omega) = \underline{a}^2(\omega)\underline{\delta}(\omega)$$

*and*

$$\overline{\overline{\beta}}(\omega) = \max\left\{\overline{a}^2(\omega)\left(1 + N\bar{\delta}(\omega)\right), \overline{a}^2(\omega)N\bar{\delta}(\omega)\left(1 + N\bar{\delta}(\omega)\right)\right\},$$

$$\overline{a}(\omega) = \exp\left(\int_0^2 \sum_{l=1}^N \max_{1 \leq k \leq N} a_{lk}(\theta_\tau\omega)\, d\tau\right), \quad \bar{\delta}(\omega) = \max_{\substack{0 \leq t \leq 2 \\ i,j=1,2,\ldots,N}} b_{ij}(\theta_t\omega).$$

(iii) $\ln^- \underline{\underline{\beta}} \in L_1((\Omega, \mathfrak{F}, \mathbb{P}))$, where $\underline{\underline{\beta}}$ is as in (ii).

For an analog of (OA3) and (OA4) for quasi-periodic systems of delay differential equations, see [21].

In the rest of this section, $\Phi$ denotes $\left((U_\omega(t))_{\omega \in \Omega, t \in [0,\infty)}, (\theta_t)_{t \in \mathbb{R}}\right)$, the measurable linear skew-product semiflow generated by (4.1) on $C([-1, 0], \mathbb{R}^N)$ covering $(\theta_t)$. The following are the main theorems of this section.

**Theorem 4.1** (Entire positive solution) *Assume* (OA1) *and* (OA2). *Assume, moreover, that $B(\omega)$ is nonsingular for each $\omega \in \Omega$ and*

$$\lim_{t\to\infty} \frac{\ln \|U_\omega(t)\|}{t} > -\infty \quad \mathbb{P}\text{-a.e. on } \Omega. \tag{4.6}$$

*Then for $\mathbb{P}$-a.e. $\omega \in \Omega$ there exists a nontrivial entire positive solution of $(4.1)_\omega$.*

In view of Proposition 2.2, a sufficient condition for the satisfaction of (4.6) is given in Theorem 4.2(2) below.

**Theorem 4.2** (1) (Generalized principal Floquet subspaces and Lyapunov exponent) *Let* (OA1) *and* (OA2) *be satisfied. Moreover, assume* (OA3)(i)–(ii) *or* (OA4)(i)–(ii). *Then $\Phi$ admits families of generalized principal Floquet subspaces $\{\tilde{E}_1(\omega)\}_{\omega \in \tilde{\Omega}_1} = \{\mathrm{span}\,\{w(\omega)\}\}_{\omega \in \tilde{\Omega}_1}$.*

(2) (Finiteness of generalized principal Lyapunov exponent) *Let* (OA1) *and* (OA2) *be satisfied. Moreover, assume* (OA3)(i)–(iii) *or* (OA4) (i)–(iii). *Then $\tilde{\lambda}_1 > -\infty$, where $\tilde{\lambda}_1$ is the generalized principal Lyapunov exponent of $\Phi$ associated to the generalized principal Floquet subspaces $\{\tilde{E}_1(\omega)\}_{\omega \in \tilde{\Omega}_1}$.*

To prove the above theorems, we first prove some propositions.

**Proposition 4.3** (Positivity) *Assume* (OA1). *Then $\Phi$ satisfies* (C2).

*Proof* It suffices to prove that $u(t; \omega, u_0) \geq 0$ for all $t \geq 0$ provided that $u_0 \in C([-1, 0], \mathbb{R}^N)^+$. Since $U_\omega^{(0)}(t)$ takes, for $t \geq 0$, $(\mathbb{R}^N)^+$ into $(\mathbb{R}^N)^+$ (see, e.g., [28, Theorem 1]), we have, by (4.3), that $u(t; \omega, u_0) \geq 0$ for all $t \in [0, 1]$. Now we proceed by induction: If $u(t; \omega, u_0) \geq 0$ for $t \in [0, k]$ for some $k \in \mathbb{N}$, repeating the previous reasoning with $\omega$ replaced by $\theta_k\omega$ we obtain that $u(t; \omega, u_0) \geq 0$ for $t \in [k, k + 1]$. $\square$

**Proposition 4.4** (Integrability). *Assume* (OA1) *and* (OA2). *Then for any* $k \in \mathbb{N}$, (C1)(i) *with* 1 *replaced by* $k$ *is satisfied for* $\Phi$, *that is, the functions*

$$\left[ \Omega \ni \omega \mapsto \sup_{0 \le s \le k} \ln^+ \|U_\omega(s)\| \in [0, \infty) \right] \in L_1((\Omega, \mathfrak{F}, \mathbb{P}))$$

*and*

$$\left[ \Omega \ni \omega \mapsto \sup_{0 \le s \le k} \ln^+ \|U_{\theta_s \omega}(k - s)\| \in [0, \infty) \right] \in L_1((\Omega, \mathfrak{F}, \mathbb{P})).$$

*Proof* Fix $\omega \in \Omega$ and $u_0 \in X^+$ with $\|u_0\| = 1$, and denote $u(\cdot) = (u_1(\cdot), \dots, u_N(\cdot)) := u(\cdot; \omega, u_0)$. For each $1 \le i \le N$ we estimate

$$\frac{du_i(t)}{dt} = \sum_{j=1}^{N} a_{ij}(\theta_t \omega) u_j(t) + \sum_{j=1}^{N} b_{ij}(\theta_t \omega) u_j(t - 1)$$

$$\le \max_{1 \le j \le N} a_{ij}(\theta_t \omega) \cdot \sum_{k=1}^{N} u_k(t) + \max_{1 \le j \le N} b_{ij}(\theta_t \omega) \cdot \sum_{k=1}^{N} u_k(t - 1).$$

Consequently, in view of Proposition 4.3,

$$\frac{d}{dt} \|u(t)\|_1 \le \sum_{i=1}^{N} \max_{1 \le j \le N} a_{ij}(\theta_t \omega) \cdot \|u(t)\|_1 + \sum_{i=1}^{N} \max_{1 \le j \le N} b_{ij}(\theta_t \omega) \cdot \|u(t - 1)\|_1,$$

for all $t \ge 0$, which implies that for any $t \ge 0$,

$$\|u(t)\|_1 \le \exp\left( \int_0^t \sum_{i=1}^{N} \max_{1 \le j \le N} a_{ij}(\theta_s \omega) \, ds \right) \|u(0)\|_1$$

$$+ \int_0^t \exp\left( \int_\tau^t \sum_{i=1}^{N} \max_{1 \le j \le N} a_{ij}(\theta_s \omega) \, ds \right) \left( \sum_{i=1}^{N} \max_{1 \le j \le N} b_{ij}(\theta_\tau \omega) \right)$$

$$\times \|u(\tau - 1)\|_1 \, d\tau. \tag{4.7}$$

By (4.7) and (OA1), for $0 \le t \le 1$,

$$\|u(t)\|_1$$

$$\le \exp\left( \int_0^1 \sum_{i=1}^{N} \max_{1 \le j \le N} a_{ij}(\theta_s \omega) \, ds \right) \cdot \left( 1 + \int_0^1 \sum_{i=1}^{N} \max_{1 \le j \le N} b_{ij}(\theta_\tau \omega) \, d\tau \right)$$

$$\max_{-1 \le \tau \le 0} \|u(\tau)\|_1 \tag{4.8}$$

which implies that for $0 \le t \le 1$,

$$\|U_\omega(t) u_0\|_1 \le \exp\left( \int_0^1 \sum_{i=1}^{N} \max_{1 \le j \le N} a_{ij}(\theta_s \omega) \, ds \right) \cdot \left( 1 + \int_0^1 \sum_{i=1}^{N} \max_{1 \le j \le N} b_{ij}(\theta_\tau \omega) \, d\tau \right) \cdot \|u_0\|_1. \tag{4.9}$$

As $C([-1, 0], \mathbb{R}^N)$ is a Banach lattice, with $\|\cdot\|_1$ a lattice norm, any $u_0 \in X$, can be written as $u_0^+ - u_0^-$ with $\|u_0^+\|_1 \leq \|u_0\|_1$, $\|u_0^-\|_1 \leq \|u_0\|_1$. By (4.9), we have that for $0 \leq t \leq 1$,

$$\|U_\omega(t)u_0\| \leq \|U_\omega(t)u_0\|_1 \leq \|U_\omega(t)u_0^+\|_1 + \|U_\omega(t)u_0^-\|_1$$

$$\leq 2\exp\left(\int_0^1 \sum_{i=1}^N \max_{1\leq j \leq N} a_{ij}(\theta_s\omega)\,ds\right) \cdot \left(1 + \int_0^1 \sum_{i=1}^N \max_{1\leq j\leq N} b_{ij}(\theta_\tau\omega)\,d\tau\right)$$

$$\times \|u_0\|_1$$

$$\leq 2\sqrt{N}\exp\left(\int_0^1 \sum_{i=1}^N \max_{1\leq j \leq N} a_{ij}(\theta_s\omega)\,ds\right) \cdot \left(1 + \int_0^1 \sum_{i=1}^N \max_{1\leq j\leq N} b_{ij}(\theta_\tau\omega)\,d\tau\right)$$

$$\times \|u_0\| \tag{4.10}$$

for all $\omega \in \Omega$, $u_0 \in X$ and $0 \leq t \leq 1$. This implies that

$$\ln^+ \|U_\omega(t)\|$$

$$\leq \ln 2 + \tfrac{1}{2}\ln N + \int_0^1 \left(\sum_{i=1}^N \max_{1\leq j\leq N} a_{ij}(\theta_\tau\omega)\right)d\tau + \ln\left(1 + \int_0^1 \sum_{i=1}^N \max_{1\leq j\leq N} b_{ij}(\theta_\tau\omega)\,d\tau\right) \tag{4.11}$$

for all $\omega \in \Omega$ and $0 \leq t \leq 1$. Observe that for $0 \leq t \leq 1$,

$$\|U_\omega(t+k-1)\| = \|U_{\theta_{k-1}\omega}(t)U_{\theta_{k-2}\omega}(1)U_{\theta_{k-2}\omega}(1)\ldots U_\omega(1)\|$$

$$\leq \|U_{\theta_{k-1}\omega}(t)\| \cdot \|U_{\theta_{k-2}\omega}(1)\| \cdot \ldots \cdot \|U_\omega(1)\|. \tag{4.12}$$

Equations (4.10), (4.11) and (4.12) together with (OA2) imply the first statement in (1). The second statement is proved in much the same way.                                                              □

**Proposition 4.5** (Focusing)

(1) *Assume* (OA1) *and* (OA3)(i)–(ii)*. Then* (C3) *with* 1 *replaced by* $N + 2$ *holds.*
(2) *Assume* (OA1) *and* (OA4)(i)–(ii)*. Then* (C3) *with* 1 *replaced by* 2 *holds.*

*Proof* (1) Assume (OA1) and (OA3)(i)–(ii). We let **e** denote both the vector $(1, 1, \ldots, 1)$ and the function constantly equal to that vector on $[-1, 0]$. Let $u_0 \in C([-1, 0], \mathbb{R}^N)^+$ be given.

We first claim that

$$\underline{\beta}(\omega)\left(\|u_0(0)\|_1 + \int_0^1 \|u_0(\tau-1)\|_1\,d\tau\right)\mathbf{e} \leq U_\omega(N+2)u_0. \tag{4.13}$$

Observe that it suffices to prove (4.13) for the case that $u_0(t) = (0, \ldots, 0, u_{0i}(t), 0, \ldots, 0)$. Without loss of generality, we assume that $u_0(t) = (u_{01}(t), 0, \ldots, 0)$. By (OA3)(i),(ii), there are $j_0, j_1, \ldots, j_{N-1}$ such that

$$b_{j_{l+1}j_l}(\omega) \geq \underline{\delta}(\omega) \quad \forall j = 0, 1, \ldots, N-1,$$

where $j_0 = 1$, $j_l \neq j_k$ for $k \neq l$, $l, k = 1, 2, \ldots, N-1$. We show by induction that for $k \leq t \leq k+1$, $k = 1, 2, \ldots, N-1$,

$$
\begin{cases}
u_{j_i}(t; \omega, u_0) \geq \underline{a}^{k+1}(\omega)\underline{\delta}^i(\omega) \left( \int_0^1 u_{01}(\tau - 1)\, d\tau \cdot \frac{1}{(i-1)!} + u_{01}(0) \cdot \frac{1}{i!} \right), \quad 1 \leq i \leq k-1 \\
u_{j_k}(t; \omega, u_0) \geq \underline{a}^{k+1}(\omega)\underline{\delta}^k(\omega) \left( \int_0^1 u_{01}(\tau - 1)\, d\tau \cdot \frac{(t-k)^{k-1}}{(k-1)!} + u_{01}(0) \cdot \frac{(t-k)^k}{k!} \right).
\end{cases}
\tag{4.14}
$$

Put $u_i(t) = u_i(t; \omega, u_0)$. Note that for $t \geq 0$,

$$
\dot{u}_1(t) \geq a_{11}(\theta_t \omega) u_1(t) \tag{4.15}
$$

and

$$
\dot{u}_{j_l}(t) \geq a_{j_l j_l}(\theta_t \omega) u_{j_l}(t) + b_{j_l j_{l-1}}(\theta_t \omega) u_{j_l}(t-1) \quad \forall l = 1, 2, \ldots, N. \tag{4.16}
$$

Hence, for $0 \leq t \leq 1$,

$$
u_1(t) \geq \exp\left( \int_0^t a_{11}(\theta_\tau \omega)\, d\tau \right) u_{01}(0) \geq \underline{a}(\omega) u_{01}(0) \tag{4.17}
$$

and

$$
u_{j_1}(t) \geq \int_0^t \exp\left( \int_\tau^t a_{j_1 j_1}(\theta_s \omega)\, ds \right) b_{j_1 1}(\theta_\tau \omega) u_1(\tau - 1)\, d\tau \geq \underline{a}(\omega)\underline{\delta}(\omega) \int_0^t u_{01}(\tau - 1)\, d\tau. \tag{4.18}
$$

By (4.15)–(4.18), for $1 \leq t \leq 2$,

$$
\begin{aligned}
u_{j_1}(t) &\geq \exp\left( \int_1^t a_{j_1 j_1}(\theta_\tau \omega)\, d\tau \right) u_{j_1}(1) + \int_1^t \exp\left( \int_\tau^t a_{j_1 j_1}(\theta_s \omega)\, ds \right) \\
&\quad \times b_{j_1 1}(\theta_\tau \omega) u_1(\tau - 1)\, d\tau \\
&\geq \underline{a}^2(\omega)\underline{\delta}(\omega) \int_0^1 u_{01}(\tau - 1)\, d\tau + \underline{a}(\omega)\underline{\delta}(\omega) \int_1^t u_1(\tau - 1)\, d\tau \\
&\geq \underline{a}^2(\omega)\underline{\delta}(\omega) \left( \int_0^1 u_{01}(\tau - 1)\, d\tau + u_{01}(0) \cdot (t-1) \right).
\end{aligned}
\tag{4.19}
$$

Hence (4.14) holds for $k = 1$.

Assume that (4.14) holds for $k = l - 1$. Then by (4.16), for $l \le t \le l + 1$,

$$u_{j_i}(t) \ge \exp\left(\int\limits_l^t a_{j_i j_i}(\theta_\tau \omega) \, d\tau\right) u_{j_i}(l)$$

$$\ge \underline{a}(\omega) u_{j_i}(l)$$

$$\ge \underline{a}^{l+1}(\omega) \underline{\delta}^i(\omega) \left(\int\limits_0^1 u_{01}(\tau - 1) \, d\tau \cdot \frac{1}{(i-1)!} + u_{01}(0) \cdot \frac{1}{i!}\right)$$

for $i = 1, 2, \ldots, l - 1$,

$$u_{j_l}(t) \ge \int\limits_l^t \exp\left(\int_\tau^t a_{j_l j_l}(\theta_s \omega) \, ds\right) b_{j_l j_{l-1}}(\theta_\tau \omega) u_{j_{l-1}}(\tau - 1) \, d\tau$$

$$\ge \underline{a}(\omega) \underline{\delta}(\omega) \int\limits_l^t u_{j_{l-1}}(\tau - 1) \, d\tau$$

$$\ge \underline{a}^{l+1}(\omega) \underline{\delta}^l(\omega) \left(\int\limits_0^1 u_{01}(\tau - 1) \, d\tau \cdot \int\limits_l^t \frac{(s-l)^{l-2}}{(l-2)!} \, ds + u_{01}(0) \int\limits_l^t \frac{(s-l)^{k-1}}{(l-1)!} \, ds\right)$$

$$= \underline{a}^{l+1}(\omega) \underline{\delta}^l(\omega) \left(\int\limits_0^1 u_{01}(\tau - 1) \, d\tau \cdot \frac{(t-l)^{l-1}}{(l-1)!} + u_{01}(0) \cdot \frac{(t-l)^l}{l!}\right).$$

By induction, (4.14) holds for $k = 1, 2, \ldots, N - 1$.

Putting in (4.14) $k = N - 1$ and $t = N$ we have

$$u_{j_i}(N) \ge \underline{a}^N(\omega) \underline{\delta}^i(\omega) \left(\int\limits_0^1 u_{01}(\tau - 1) \, d\tau \cdot \frac{1}{(i-1)!} + u_{01}(0) \cdot \frac{1}{i!}\right), \quad 1 \le i \le N - 1,$$
$$(4.20)$$

from which it follows, noting that $\dot{u}_{j_i}(t) \ge a_{j_i j_i}(\theta_t \omega) u_{j_i}(t)$, that

$$u_{j_i}(t) \ge \underline{a}^{N+2}(\omega) \underline{\delta}^i(\omega) \left(\int\limits_0^1 u_{01}(\tau - 1) \, d\tau \cdot \frac{1}{(i-1)!} + u_{01}(0) \cdot \frac{1}{i!}\right),$$
$$N \le t \le N + 2, \quad 1 \le i \le N - 1. \tag{4.21}$$

Regarding $u_1(t)$, observe that there is $k \in \{1, \ldots, N-1\}$ such that $b_{1j_k}(\theta_t\omega) \geq \underline{\delta}(\omega)$ for $0 \leq t \leq N+2$. For $N \leq t \leq N+2$, by (4.21), we have

$$
\begin{aligned}
u_1(t) &\geq \int_N^t \exp\left(\int_\tau^t a_{11}(\theta_s\omega)\, ds\right) b_{1j_k}(\theta_\tau\omega) u_{j_k}(\tau-1)\, d\tau \\
&\geq \underline{a}(\omega)\underline{\delta}(\omega) \int_N^t u_{j_k}(\tau-1)\, d\tau \\
&\geq \underline{a}^{N+3}(\omega)\underline{\delta}^{k+1}(\omega)\left(\int_0^1 u_{01}(\tau-1)\, d\tau \cdot \frac{1}{(k-1)!} + u_{01}(0) \cdot \frac{1}{k!}\right)(t-N). \quad (4.22)
\end{aligned}
$$

It follows from (4.22) that

$$
u_1(t) \geq \underline{a}^{N+3}(\omega)\underline{\delta}^{k+1}(\omega)\left(\int_0^1 u_{01}(\tau-1)\, d\tau \cdot \frac{1}{(k-1)!} + u_{01}(0) \cdot \frac{1}{k!}\right), \quad N+1 \leq t \leq N+2. \tag{4.23}
$$

Putting together (4.21) and (4.23) we obtain (4.13).

Next, we claim that for any $u_0 \in C([-1,0], \mathbb{R}^N)^+$,

$$
u(t; \omega, u_0) \leq \bar{\beta}(\omega)\left(\int_0^1 \|u_0(\tau-1)\|_1\, d\tau + \|u_0(0)\|_1\right)\mathbf{e}, \quad N+1 \leq t \leq N+2. \tag{4.24}
$$

To prove (4.24), let $v(t) = \|u(t; \omega, u_0)\|_1$. Notice that for $t \geq 0$, we have

$$
\dot{v}(t) \leq \sum_{l=1}^N \max_{1 \leq k \leq N} a_{lk}(\theta_t\omega) \cdot v(t) + N\bar{\delta}(\omega)v(t-1). \tag{4.25}
$$

We prove by induction that for $k-1 \leq t \leq k$, $k = 1, 2, \ldots, N+2$,

$$
v(t) \leq \bar{a}^k(\omega)(1 + N\bar{\delta}(\omega))^{k-1}v(0) + \bar{a}^k(\omega)N\bar{\delta}(\omega)(1 + N\bar{\delta}(\omega))^{k-1}\int_0^1 v(\tau-1)\, d\tau. \tag{4.26}
$$

In fact, for $0 \leq t \leq 1$, by (4.25),

$$
\begin{aligned}
v(t) &\leq \exp\left(\int_0^t \sum_{l=1}^N \max_{1 \leq k \leq N} a_{lk}(\theta_\tau\omega)\, d\tau\right) v(0) \\
&\quad + N\bar{\delta}(\omega)\int_0^t \exp\left(\int_\tau^t \sum_{l=1}^N \max_{1 \leq k \leq N} a_{lk}(\theta_s\omega)\, ds\right) v(\tau-1)\, d\tau \\
&\leq \bar{a}(\omega)v(0) + \bar{a}(\omega)N\bar{\delta}(\omega)\int_0^1 v(\tau-1)\, d\tau.
\end{aligned}
$$

Hence (4.26) holds for $k = 1$. Assume that (4.26) holds for $k = l$. Then for $l \leq t \leq l + 1$,

$$
\begin{aligned}
v(t) &\leq \exp\left(\int_l^t \sum_{l=1}^N \max_{1 \leq k \leq N} a_{lk}(\theta_\tau \omega)\, d\tau\right) v(l) \\
&\quad + N\bar{\delta}(\omega) \int_l^t \exp\left(\int_\tau^t \sum_{l=1}^N \max_{1 \leq k \leq N} a_{lk}(\theta_s \omega)\, ds\right) v(\tau - 1)\, d\tau \\
&\leq \bar{a}(\omega) v(l) + \bar{a}(\omega) N\bar{\delta}(\omega) \int_l^t v(\tau - 1)\, d\tau \\
&\leq \bar{a}(\omega)(1 + N\bar{\delta}(\omega)) \Big[ \bar{a}^l(\omega)(1 + N\bar{\delta}(\omega))^{l-1} v(0) \\
&\quad + \bar{a}^l(\omega) N\bar{\delta}(\omega)(1 + N\bar{\delta}(\omega))^{l-1} \int_0^1 v(\tau - 1)\, d\tau \Big] \\
&= \bar{a}^{l+1}(\omega)(1 + N\bar{\delta}(\omega))^l v(0) + \bar{a}^{l+1}(\omega) N\bar{\delta}(\omega)(1 + N\bar{\delta}(\omega))^l \int_0^1 v(\tau - 1)\, d\tau.
\end{aligned}
$$

By induction, (4.26) holds. The claim (4.24) follows from (4.26). (1) follows from (4.13) and (4.24).

(2) Assume (OA1) and (OA4)(i)–(ii).

We first prove that for any $u_0 \in C([-1, 0], \mathbb{R}^N)^+$,

$$
U_\omega(2) u_0 \geq \underline{\underline{\beta}}(\omega) \left( \|u_0(0)\|_1 + \int_0^1 \|u_0(\tau - 1)\|_1\, d\tau \right) \mathbf{e}. \tag{4.27}
$$

Note that

$$
\dot{u}_i(t; \omega, u_0) \geq a_{ii}(\theta_t \omega) u_i(t; \omega, u_0) + \underline{\delta}(\omega) \|u(t - 1; \omega, u_0)\|_1 \quad \forall i = 1, 2, \ldots, N.
$$

Hence for $0 \leq t \leq 1$ and $1 \leq i \leq N$,

$$
\begin{aligned}
u_i(t; \omega, u_0) &\geq \exp\left(\int_0^t a_{ii}(\theta_\tau \omega)\, d\tau\right) u_{0i}(0) + \int_0^t \exp\left(\int_\tau^t a_{ii}(\theta_s \omega)\, ds\right) \\
&\quad \times \underline{\delta}(\omega) \|u_0(\tau - 1)\|_1\, d\tau \\
&\geq \underline{a}(\omega) u_{0i}(0) + \underline{a}(\omega) \underline{\delta}(\omega) \int_0^t \|u_0(\tau - 1)\|_1\, d\tau.
\end{aligned}
$$

Then for $1 \leq i \leq N$ and $1 \leq t \leq 2$, we have

$$
\begin{aligned}
u_i(t; \omega, u_0) &\geq \exp\left(\int_1^t a_{ii}(\theta_\tau \omega)\, d\tau\right) u_i(1; \omega, u_0) \\
&\quad + \int_1^t \exp\left(\int_\tau^t a_{ii}(\theta_s \omega)\, ds\right) \underline{\delta}(\omega) \|u(\tau - 1; \omega, u_0)\|_1\, d\tau \\
&\geq \underline{a}(\omega)\left(\underline{a}(\omega)\|u_0(0)\|_1 + \underline{a}(\omega)\underline{\delta}(\omega)\int_0^1 \|u_0(\tau - 1)\|_1\, d\tau\right) \\
&\quad + \underline{a}(\omega)\underline{\delta}(\omega)\int_1^t \|u(\tau - 1; \omega, u_0)\|_1\, d\tau \\
&\geq \underline{a}^2(\omega)\|u_0(0)\|_1 + \underline{a}^2(\omega)\underline{\delta}(\omega)\int_0^1 \|u_0(\tau - 1)\|_1\, d\tau.
\end{aligned}
$$

Hence (4.27) holds.

Next, we note that, by the arguments as in (1), there holds

$$
u(t; \omega, u_0) \leq \bar{\bar{\beta}}(\omega)\left(\|u_0(0)\|_1 + \int_0^1 \|u_0(\tau - 1)\|_1\, d\tau\right)\mathbf{e}, \quad 1 \leq t \leq 2. \tag{4.28}
$$

(2) follows from (4.27) and (4.28). □

*Proof of Theorem 4.1* By Proposition 4.2, $U_\omega(k)$ satisfies (C1)(iii). Assume $B(\omega)$ is non-singular. Then $U_\omega(k)$ is injective and hence (C1)(ii) is satisfied. By Proposition 4.4, (C1)(i) is satisfied with 1 replaced by any $k \in \mathbb{N}$. By Proposition 4.3, (C2) is satisfied with 1 replaced by any $k \in \mathbb{N}$. Hence the conditions in Proposition 2.3 are satisfied. Theorem 4.1 then follows from Theorem 2.1. □

*Proof of Theorem 4.2* (1) In view of Propositions 4.4 through 4.5 the first part follows from Theorem 2.2(5).

(2) Assume, for definiteness, that (OA3)(iii) is satisfied. It follows from the proof of Proposition 4.5 (Eq. (4.13)) and Proposition 4.3 that

$$
U_\omega(N + 2)\mathbf{e} \geq 2N\underline{\beta}(\omega)\mathbf{e},
$$

consequently, by the monotonicity of the norm $\|\cdot\|$,

$$
\frac{\ln \|U_\omega(n(N+2))\mathbf{e}\|}{n(N+2)} \geq \frac{1}{n(N+2)}\sum_{k=0}^{n-1} \ln \underline{\beta}(\theta_{k(N+2)}\omega) + \frac{\ln(2N)}{N+2}.
$$

An application of the Birkhoff ergodic theorem to the function $-\ln \underline{\beta}$ gives that for $\mathbb{P}$-a.e. $\omega \in \Omega$ there holds

$$
\lim_{n \to \infty} \frac{1}{n}\sum_{k=0}^{n-1} \ln \underline{\beta}(\theta_{k(N+2)}\omega)) = \int_\Omega \ln \underline{\beta}(\theta_{N+2}\cdot)\, d\mathbb{P}(\cdot) > -\infty,
$$

hence

$$\limsup_{t \to \infty} \frac{\ln \|U_\omega(t)\mathbf{e}\|}{t} > -\infty.$$

By Definition 2.2(iv), for $\mathbb{P}$-a.e. $\omega \in \Omega$ we have

$$\limsup_{t \to \infty} \frac{\ln \|U_\omega(t)\mathbf{e}\|}{t} \le \tilde{\lambda}_1,$$

which concludes our proof.                                                □

In systems of delay differential equations the choice of $C([-1, 0], \mathbb{R}^N)$ as the "state space" is not the only possible: observe that, since the dual space is not separable, we are unable to apply the theory of generalized exponential separation as presented in [18]. It appears that applying the (separable and reflexive) space $L_2((-1, 0), \mathbb{R}^N) \oplus \mathbb{R}^N$ (as in [3]) could be useful in proving such properties.

For some linear time-periodic (systems of) delay differential equations with an additional structure invariant decompositions into countably many finite-dimensional subbundles (labelled by a lap number) were proved in [14,15].

**Acknowledgments**    Janusz Mierczyński is supported by the NCN Grant Maestro 2013/08/A/ST1/00275.